\documentclass[11pt,a4paper]{article}
\usepackage[OT1]{fontenc}
\usepackage[english]{babel}
\usepackage[round]{natbib}
\let\cite\citet

\usepackage{authblk}
\usepackage{amsthm,amsmath,amsfonts,amssymb}
\usepackage[colorlinks,linkcolor=blue,citecolor=blue,urlcolor=blue,]{hyperref}
\usepackage{graphicx,color}
\usepackage[dvipsnames]{xcolor}
\usepackage[normalem]{ulem}
\usepackage{subcaption}
\usepackage{mathrsfs}
\usepackage{bm,dsfont,ifsym}
\usepackage{xhfill}

\allowdisplaybreaks  
\theoremstyle{plain}

\newtheorem{theorem}{Theorem}[section]

\newtheorem{proposition}[theorem]{Proposition}

\theoremstyle{remark}
\newtheorem{definition}[theorem]{Definition}
\newtheorem*{example}{Example}

\newtheorem*{rmk}{Remark}

\def\eqd{\stackrel{\mbox{\scriptsize{d}}}{=}}

\definecolor{mycolor1}{rgb}{0.6350, 0.0780, 0.1840}
\definecolor{mycolor2}{rgb}{0, 0.4470, 0.7410}
\definecolor{mycolor3}{rgb}{0.9290, 0.6940, 0.1250}
\definecolor{mycolor4}{rgb}{0.4660, 0.6740, 0.1880}
\definecolor{mycolor5}{rgb}{0.8500, 0.3250, 0.0980}

\begin{document}
	\begin{center}
		\textmd{\Large{\bfseries{{Ferguson's Dirichlet Process Breakthrough: A Lasting Legacy}}}}
	\end{center}
	\begin{center}
		{A. Lijoi$^1$, I. Pr\"unster$^1$ and J. Zhang$^2$}
	\end{center}

	\begin{quote}
		\begin{small}
			\begin{center}
				\noindent $^1$ Bocconi Institute for Data Science and Analytics, Bocconi University, Milan, Italy.\\ \textit{E-mail}: antonio.lijoi@unibocconi.it; igor@unibocconi.it
				
				\noindent $^2$ Department of Mathematics and Information Technology, The Education University of Hong Kong, Hong Kong.\\ \textit{E-mail}: zjunyi@eduhk.hk
			\end{center}
		\end{small}
	\end{quote}
	
	\smallskip
	
	\begin{abstract}
		Ferguson's 1973 introduction of the Dirichlet process marked a breakthrough in Bayesian nonparametric statistics. For the first time, a prior on the space of probability measures fulfilled two key \textit{desiderata}: large support and analytical tractability. In this paper, we review three complementary constructions of the Dirichlet process, whose roots can be traced back to Ferguson: through finite-dimensional distributions, via normalization of a gamma process, and through predictive distributions. Each perspective not only deepens the understanding of the Dirichlet process but also provides a template for generalizations, from normalized random measures with independent increments to Gibbs--type priors and beyond. Over the past fifty years, the Dirichlet process has become the cornerstone of Bayesian nonparametric methodology and applications, while simultaneously inspiring the expansion of the landscape of nonparametric priors. Since de~Finetti laid out the Bayesian nonparametric framework in the 1930s, the key obstacle had been the absence of a tractable nonparametric prior. Ferguson's contribution overcame this challenge, providing a solution to a decades-long open problem. In recognition of this decisive advance, it seems appropriate to refer to the Dirichlet process as the Ferguson--Dirichlet process.
		
		\par
		
		\vspace{9pt}
		\par
	\end{abstract}

\section{Introduction}
\subsection{Exchangeability and representation theorem}
		
Learning from experience, and induction more generally, rests on the philosophical principle that the future will resemble the past. This idea, known as the \textit{Uniformity of Nature Principle}, originates in Hume's analysis of induction \citep{hume1739} and has been discussed extensively in the philosophy of science literature; see e.g. \cite{salmon1953}. De~Finetti \citep{definetti1937} embedded this notion of induction within probability theory by reformulating the required symmetry between past and future as a probabilistic invariance.
\begin{definition}[Exchangeability, \citealt{definetti1937}]
	A sequence $(X_n)_{n \geq 1}$ of $\mathbb{X}$-valued random variables is exchangeable if for any $n\geq 1$ and permutation $\pi$ of $(1, \dots, n)$, 
	\begin{equation*}
		(X_1, \dots, X_n) \stackrel{d}{=} (X_{\pi(1)}, \dots, X_{\pi(n)}).
	\end{equation*}
\end{definition}
Exchangeability is weaker than independence, and it is precisely this feature that enables principled prediction: under strict independence no information about the future can be gained from the past, whereas exchangeability allows learning to take place. De Finetti's representation theorem provides a rigorous foundation for this principle by showing that every exchangeable sequence can be viewed as conditionally i.i.d. given a random probability measure.
\begin{theorem}[de~Finetti's representation theorem]
	A sequence of random variables $(X_n)_{n \geq 1}$ taking values in a measurable space $\mathbb{X}$ is exchangeable if and only if there exists a probability measure $Q$ on the space $\mathscr{P}$\footnote{Throughout we assume that $\mathbb{X}$, endowed with the corresponding Borel $\sigma$-field $\mathscr{X}$, is a complete and separable metric space, also known as Polish space. This, in turn, ensures that $\mathscr{P}$ on $\mathbb{X}$ equipped with the $\sigma$-field induced by the weak topology is also Polish and, crucially for Bayesian inference, that the conditional distributions are regular.} of probability measures on $\mathbb{X}$ such that, for all $n \geq 1$ and measurable sets $A_1, \dots, A_n \subseteq \mathbb{X}$,
	\begin{equation}\label{eq:definetti}
		\mathbb{P}(X_1 \in A_1, \dots, X_n \in A_n) = 
		\int_{\mathscr{P}} \prod_{i=1}^n P(A_i) \, Q(dP).
	\end{equation}
\end{theorem}
The mixing measure $Q$ is known as the de~Finetti measure of $(X_n)_{n\geq 1}$: it is the law of a random probability measure $\tilde P$ and serves as a prior distribution for Bayesian inference. More broadly, the representation theorem provides a fundamental justification of the Bayesian approach, showing that inductive learning from exchangeable data is equivalent to assigning a prior distribution to a parameter, which in the general case \eqref{eq:definetti} is the probability distribution itself. See also \cite[Ch.~7]{diaconis2018} for further discussion of its central role in justifying the Bayesian framework.
	
Equivalently, de~Finetti's representation theorem admits the hierarchical formulation commonly used in Bayesian statistics, where it naturally serves as a generative model
\begin{equation}\label{BayesianHierarhicalModel}
	\begin{split}
		X_i \mid \tilde{P} \: &\stackrel{\mbox{\scriptsize iid}}{\sim} \: \tilde{P} ,\quad i=1, \dots, n ,\\
		\tilde{P} \:&\sim\: Q .
	\end{split}
\end{equation}
The structure of the prior distribution $Q$ determines whether the model \eqref{BayesianHierarhicalModel} is parametric or nonparametric. The parametric case is recovered as a special instance when $Q$ is degenerate on a subset of $\mathscr{P}$ indexed by a finite-dimensional parameter, so that the prior effectively reduces to one on a conventional parameter space; for example, if $Q$ is concentrated on Gaussian distributions, the model reduces to placing a prior on the mean and the variance, $(\mu,\sigma^2)$ of the $X_i$'s. When $Q$ is not concentrated on a class of distributions indexed by a finite-dimensional parameter, the model is defined as nonparametric, thereby extending the framework to an infinite-dimensional parameter, namely a probability measure on $\mathbb{X}$. As observed by B.~de~Finetti as early as 1935 \citep{definetti1935}, in nonparametric problems
\begin{quote}
	\textsl{the curve itself is a random element that has, a priori, a certain probability distribution in the functional space, which is modified, following the observation of certain data, still on the basis of Bayes' theorem.}
\end{quote}
	
\noindent The general Bayesian framework was laid out by de~Finetti in the 1930s and  given a rigorous mathematical foundation in \cite{hewitt1955}, but the main obstacle for nonparametric problems was the absence of a tractable prior distribution on the space of probability measures. As Lindley remarked in his review of Bayesian statistics \citep[Ch.~12]{lindley72}:
\begin{quote}
	\textsl{It is perhaps worth stopping to remark that the problem is a technical one; the Bayesian method embraces non-parametric problems but cannot solve them because the requisite tool is missing.}
\end{quote}
	
\subsection{Ferguson's breakthrough}
	
In order to address the lack of a workable prior distribution for nonparametric inference, Ferguson's seminal paper \citep{MR350949} begins by formulating two desiderata that such a prior should satisfy. As noted in the previous section, the term \emph{nonparametric} refers to cases where the mixing measure $Q$ in \eqref{eq:definetti} is not concentrated on a subset of $\mathscr{P}$ characterized by a finite-dimensional parameter. The choice of such a prior, however, is not arbitrary, as on such large spaces one must remain aware of the implications of the prior choice to ensure meaningful statistical inference. Ferguson's desiderata therefore specify the key properties of a useful nonparametric prior: sufficiently large support and analytical tractability of the posterior. He states them as follows:
\begin{quote} 
	\textsl{The Bayesian approach to statistical problems, though fruitful in many ways, has been rather unsuccessful in treating nonparametric problems. This is due primarily to the difficulty in finding workable prior distributions on the parameter space, which in nonparametric problems is taken to be a set of probability distributions on a given sample space. There are two desirable properties of a prior distribution for nonparametric problems.\\
		(I) The support of the prior distribution should be large---with respect to some suitable topology on the space of probability distributions on the sample space.\\
		(II) Posterior distributions given a sample of observations from the true probability distribution should be manageable analytically.}
\end{quote}
	
Ferguson's decisive contribution was the derivation of a nonparametric prior that simultaneously satisfies both desiderata -- large support and analytical tractability. The construction is so natural and compelling that it has stood, to this day, as the cornerstone of Bayesian Nonparametrics. This is the process widely known as the \emph{Dirichlet process}, but given the magnitude of the breakthrough it represents, it is more appropriately referred to as the \emph{Ferguson--Dirichlet process}.
	
Before \cite{MR350949}, several Bayesian constructions for random distribution functions had been proposed, but none simultaneously met the above desiderata, in particular analytical tractability. More specifically, the procedure of \cite{dubins1967} yields singular continuous laws, whereas the constructions of \cite{kraft1964,krafteeden1964} produce absolutely continuous laws but lack a generally workable posterior calculus. The class of \emph{tailfree} priors, introduced by \cite{freedman1963} on a countable space and extended by \cite{fabius1964}, does offer computable posteriors in specific settings, and the Ferguson--Dirichlet process itself can be viewed as a particularly simple tailfree law. Nevertheless, taken together these strands did not provide a broadly applicable, tractable prior on the space of probability measures, and Ferguson's contribution filled precisely this gap. In subsequent work, \cite{MR438568} substantially developed tailfree processes and introduced the P\'olya tree architecture of priors, which was further advanced by \cite{lavine1992} and \cite{mauldin1992}, leading to a stream of modern Bayesian nonparametric contributions; see, e.g., \cite{li2010,li2024}.
	
More specifically, let $\alpha$ be a finite, non-null measure on $\mathbb{X}$. A Ferguson--Dirichlet process with parameter--measure $\alpha$ is a random probability measure $\tilde P$ on $\mathbb{X}$ to be denoted as $\mathscr{D}(\alpha)$ throughout. Its law $Q$ on $\mathscr{P}$ is the de~Finetti measure of an exchangeable sequence and, in the Bayesian framework, it serves as a nonparametric prior distribution. We postpone the formal constructions of $\mathscr{D}(\alpha)$ to the following section but we emphasize already here the way in which the Ferguson--Dirichlet process satisfies the two desiderata.
	
We begin with the support property (I). The natural topology to consider on the space of probability measures $\mathscr{P}$ is the weak topology, and the Ferguson--Dirichlet process prior assigns positive probability to every weak neighborhood of any $P \in \mathscr{P}$ whose support is contained in that of $\alpha$. This is stated in \cite{MR438568} as 
\begin{quote}
	\textsl{The support of $\mathscr{D}(\alpha)$ with respect to the topology of weak convergence is the set of all distributions whose support is contained in the support of $\alpha$.}
\end{quote}
This result, first formulated in \cite{MR350949,MR438568} and later made precise in \cite{dalal1980} and \cite{majumdar1992}, is typically referred to as the Ferguson--Dirichlet process enjoying the ``full support property''. It guarantees that any potential data--generating distribution is contained in its support, therefore fully satisfying the first desideratum.
	
We now turn to the desideratum (II), namely analytical tractability, which is met by the Ferguson--Dirichlet process.
\begin{theorem}[{\citealt[Theorem~1]{MR350949}}]
	\label{TheoremDirichletProcessPosterior}
	Let $\tilde P \sim \mathscr{D}(\alpha)$ on $\mathbb{X}$, and let $X_1,\dots,X_n$ be a sample of size $n$ from $\tilde P$. Then the conditional distribution of $\tilde P$ given $X_1,\dots,X_n$ is again a Ferguson--Dirichlet process, with updated parameter measure $\alpha + \sum_{i=1}^n \delta_{X_i}$, where $\delta_x$ denotes the point mass at $x$.
\end{theorem}
In other words, the posterior distribution of a Ferguson--Dirichlet process is still a Ferguson--Dirichlet process whose parameter-measure is updated by adding a point mass at each observation of the sample. This property is typically referred to as the Ferguson--Dirichlet process being \emph{conjugate}.
	
However, it is useful to clarify the notion of conjugacy in the nonparametric framework. In the parametric setting, one usually speaks of a conjugate pair between a likelihood and a prior. In the nonparametric case the situation is subtler. In fact, following \cite{MR2730661}, one can distinguish between \emph{parametric conjugacy}, where the posterior process is of the same type as the prior with updated parameters, and \emph{structural conjugacy}, where the prior is embedded in a larger class of processes and the posterior is shown to belong to that class, even if it is not of exactly the same type.
	
While the Ferguson--Dirichlet process enjoys parametric conjugacy and is largely the only nonparametric prior to have it (see, e.g., \citealt{MR2255112}), when constructing a nonparametric prior one typically aims for structural conjugacy or even weaker forms of posterior tractability. An early and important example of structural conjugacy is the class of neutral-to-the-right processes, introduced in \cite{doksum1974} (see also \citealt{MR438568}). 
	
The paper is organized as follows. In Sections~\ref{Section2}, \ref{Section3}, and \ref{Section4}, we review, from a modern perspective, the three constructions of the Ferguson--Dirichlet process that originated in Ferguson's work, namely those based on finite-dimensional distributions, normalization of completely random measures, and predictive distributions. In addition to Ferguson's seminal 1973 paper, our review also covers, to some extent, \cite{MR373022} and \cite{MR438568}. In each section we also discuss how these constructions can be extended to derive other interesting nonparametric priors. Finally, in Section~\ref{SectionRemarks} we broaden the perspective with a brief discussion of extensions beyond exchangeability and the dependent Ferguson--Dirichlet process. Henceforth, we adhere as closely as possible to the original notation used in \cite{MR350949}, with any deviations made only to maintain consistency across sections.

\section{Construction via the finite-dimensional distributions}\label{Section2}
	
In this section, we review the construction of the Ferguson--Dirichlet process through its finite-dimensional distributions and discuss how this approach can be extended to obtain other well-known nonparametric priors.
	
\subsection{The Ferguson--Dirichlet process}
	
\begin{figure}[htbp]
	\centering
	\includegraphics[width=0.5\textwidth]{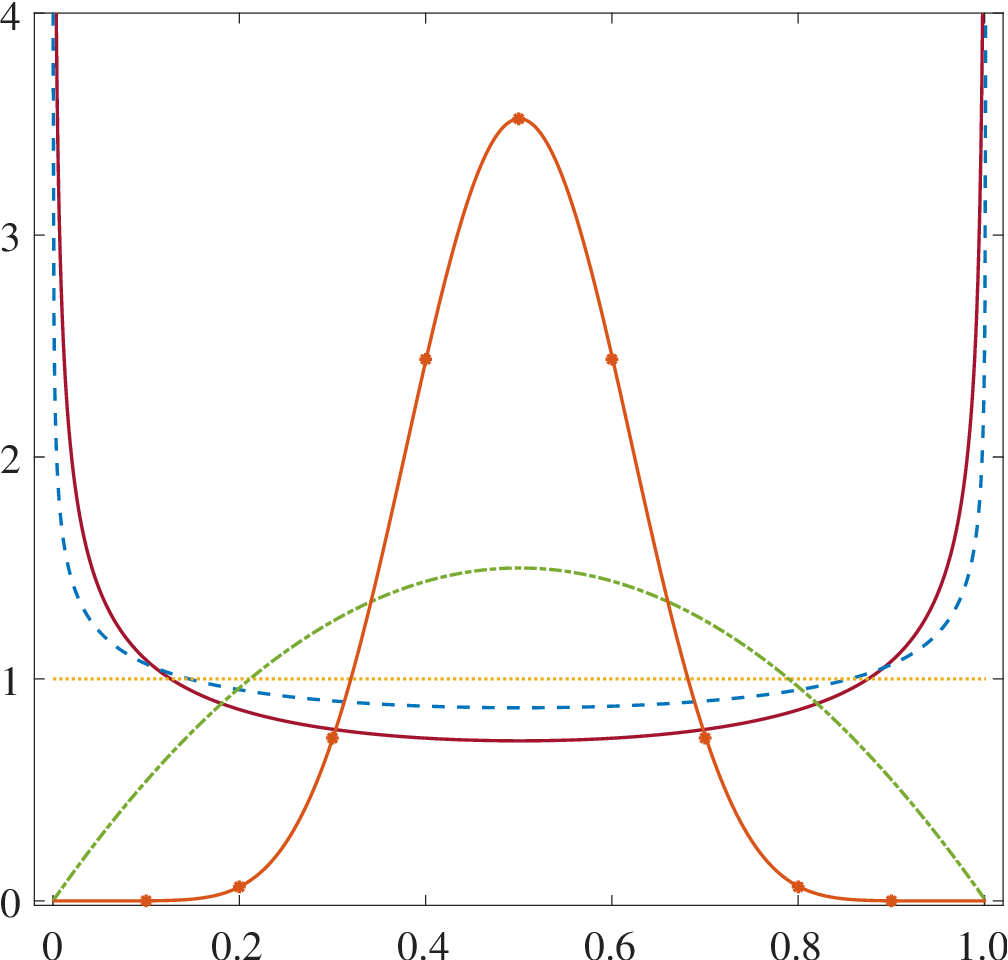}
	\caption{Beta densities on $(0,1)$ with the parameters $(\alpha_i, \alpha_i)$. \textcolor{mycolor1}{Dir$(.6, .6)$ \rule[.5ex]{1.7em}{.8pt}}; \textcolor{mycolor2}{Dir$(.8, .8)$ \rule[.5ex]{0.5em}{.8pt}\,\,\rule[.5ex]{0.5em}{.8pt}\,\,\rule[.5ex]{0.5em}{.8pt}}; \textcolor{mycolor3}{Dir$(1, 1)$ \rule[.5ex]{0.1em}{.8pt}\,\rule[.5ex]{0.1em}{.8pt}\,\rule[.5ex]{0.1em}{.8pt}\,\rule[.5ex]{0.1em}{.8pt}\,\rule[.5ex]{0.1em}{.8pt}\,\rule[.5ex]{0.1em}{.8pt}}; \textcolor{mycolor4}{Dir$(2, 2)$ \rule[.5ex]{0.5em}{.8pt}\,\rule[.5ex]{0.1em}{.8pt}\,\rule[.5ex]{0.5em}{.8pt}\,\rule[.5ex]{0.1em}{.8pt}\,\rule[.5ex]{0.5em}{.8pt}}; \textcolor{mycolor5}{Dir$(10, 10)$ \rule[.5ex]{0.8em}{.8pt}$*$\rule[.5ex]{0.8em}{.8pt}}.}
	\label{FigureDirichletDensityK2}
\end{figure}

\begin{figure}
	\centering
	\begin{subfigure}[t]{0.4\textwidth}
		\includegraphics[width=1\textwidth]{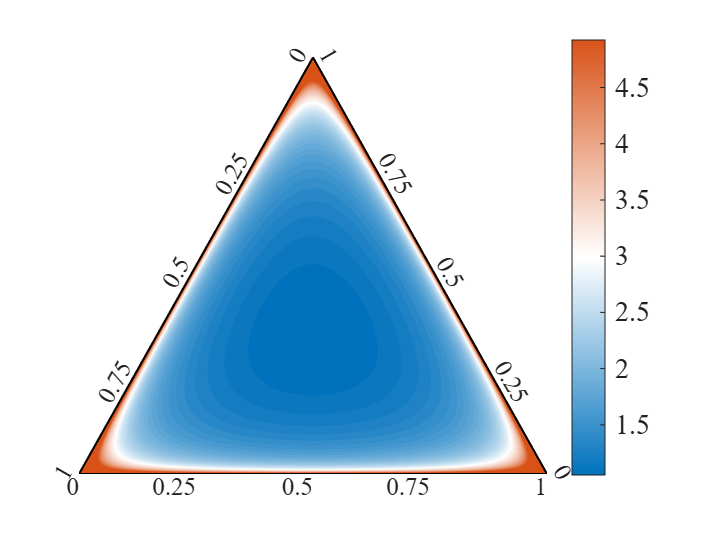}
		\caption{Dirichlet$(.6,.6,.6)$.}
	\end{subfigure}
	\enskip
	\begin{subfigure}[t]{0.4\textwidth}
		\includegraphics[width=1\textwidth]{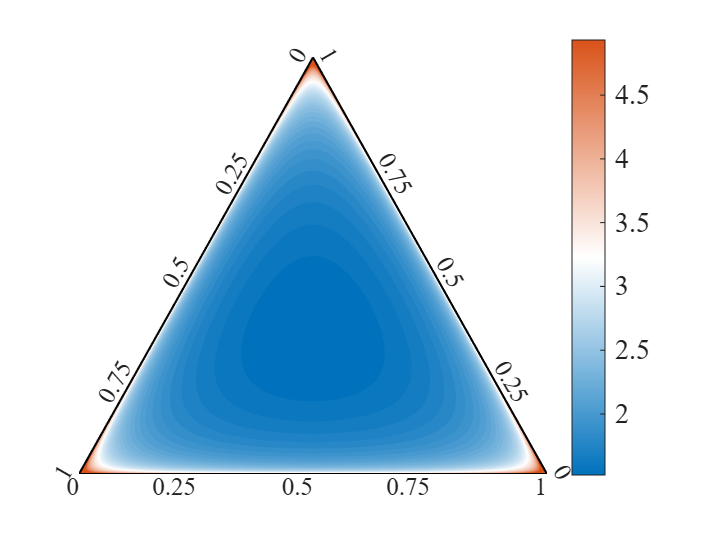}
		\caption{Dirichlet$(.8,.8,.8)$.}
	\end{subfigure}

	\begin{subfigure}[t]{0.4\textwidth}
		\includegraphics[width=1\textwidth]{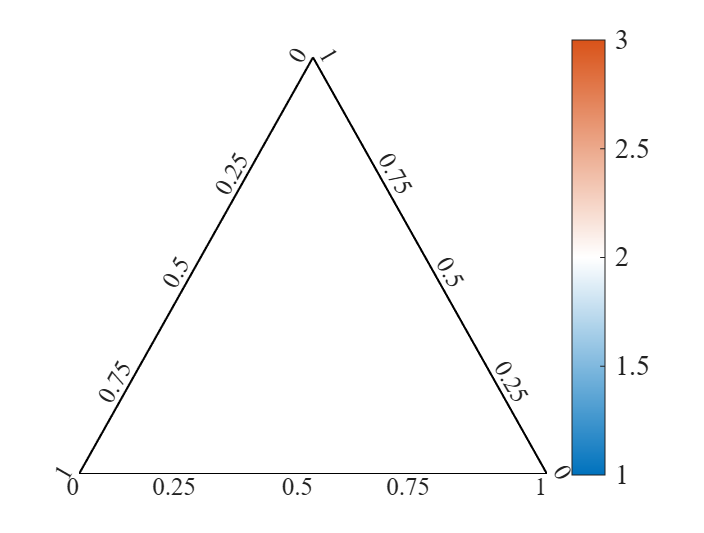}
		\caption{Dirichlet$(1,1,1)$.}
	\end{subfigure}
	\enskip
	\begin{subfigure}[t]{0.4\textwidth}
		\includegraphics[width=1\textwidth]{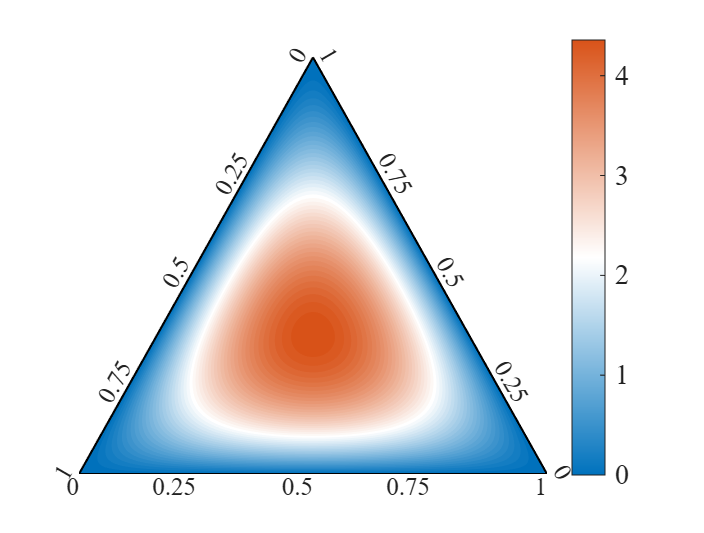}
		\caption{Dirichlet$(2,2,2)$.}
	\end{subfigure}
	\caption{Heatmaps for the Dirichlet density on $\Delta_2$ corresponding to various parameter settings.}
	\label{FigureDirichletDensityK3}
\end{figure}

We start by recalling the definition of the (finite-dimensional) Dirichlet distribution, which underlies the construction of the Ferguson--Dirichlet process. Let $\bm{\alpha}:=(\alpha_1, \ldots, \alpha_K)$ with $\alpha_i > 0$, for $i=1, \ldots, K$, be the parameters, and let $|\bm{\alpha}| := \sum_{i=1}^{K}\alpha_i$. Denote by $\Delta_{K-1}$ the $(K-1)$-dimensional simplex, namely
\begin{equation*}
	\Delta_{K-1} :=  \left\{ x_1, \ldots, x_{K-1} \geq 0 : \sum_{i=1}^{K-1}x_i \leq 1 \right\} .
\end{equation*}
The random vector $(Y_1, \dots, Y_K)$ is said to have Dirichlet distribution with parameters $\alpha_1, \dots, \alpha_K$, denoted by 
\begin{equation*}
	(Y_1, \dots, Y_K) \sim \mathscr{D}(\alpha_1, \dots, \alpha_K) ,
\end{equation*}
if 	$Y_1+\,\cdots\,+Y_K=1$ and $(Y_1,\ldots,Y_{K-1})$ has density function (with respect to the Lebesgue measure) given by
\begin{equation*}
	f_K(\bm{y}; \bm{\alpha}) = \frac{\Gamma(|\bm{\alpha}|)}{\prod_{i=1}^{K}\Gamma(\alpha_i)} \prod_{i=1}^{K-1}y_i^{\alpha_i-1} \left( 1 - |\bm{y}| \right)^{\alpha_K-1} \mathds{1}_{\Delta_{K-1}}(\bm{y}),
\end{equation*}
where $\bm{y}=(y_1,\ldots,y_{K-1})$, $|\bm{y}|:=\sum_{i=1}^{K-1} y_i$, and $\mathds{1}_A$ is the indicator function of set $A$. When $K=2$, $f_2$ becomes the density function of a Beta distribution with parameters $(\alpha_1,\alpha_2)$, namely $\mathscr{D}(\alpha_1,\alpha_2)=\mbox{beta}(\alpha_1,\alpha_2)$. Figures~\ref{FigureDirichletDensityK2} and~\ref{FigureDirichletDensityK3} illustrate the Dirichlet distribution for $(0,1)$ and $\Delta_2$ under different parameter settings.

The construction of the Ferguson--Dirichlet process in terms of finite-dimensional distributions relies on a fundamental property of the Dirichlet distribution:

\noindent \textsc{Additivity.}
Let $(Y_1, \dots, Y_K) \sim \mathscr{D}(\alpha_1, \dots, \alpha_K)$, and let $r_1, \dots, r_{l-1}$ be positive integers such that $r_1<\dots<r_{l-1}<K$. Then 
\begin{equation} \label{prop:additivity}
		\left( \sum_{i=1}^{r_1}Y_i , \sum_{i=r_1+1}^{r_2}Y_i , \dots , \sum_{i=r_{l-1}+1}^{K}Y_i \right) \sim \mathscr{D} \left( \sum_{i=1}^{r_1}\alpha_i , \sum_{i=r_1+1}^{r_2}\alpha_i , \dots , \sum_{i=r_{l-1}+1}^{K}\alpha_i \right) .
\end{equation}
\smallskip

In particular, for each $i=1,\dots, K$, the additivity property~\eqref{prop:additivity} implies that
the marginal distribution of $Y_i$ is a beta distribution, i.e.,
\begin{equation}\label{MarginalDistributionDirichletProcess}
	Y_i \sim {\text{beta}} \left( \alpha_i , \left( \sum_{j=1}^{K} \alpha_j \right) - \alpha_i \right) .
\end{equation}

In order to leverage the Dirichlet distribution as the finite-dimensional distributions of a nonparametric prior on the space of probability measures $\mathscr{P}$, Ferguson extended the definition of the Dirichlet distribution \citep{wilks62} to allow some parameters to be zero (correspondingly decreasing the dimension of the simplex on which its density is defined). He also ingeniously introduced a measure $\alpha$ that assigns the corresponding vector of parameters to each finite-dimensional distribution. This led to the following definition of the Ferguson--Dirichlet process, which we reproduce verbatim from \cite{MR350949} with minimal edits for notational consistency. 

\begin{definition}[Dirichlet process, \citealt{MR350949}]\label{DefinitionDirichletProcess1}
Let $\alpha$ be a non-null finite measure (non-negative and finitely additive) on $(\mathbb{X},\mathscr{X})$. We say $\tilde{P}$ is a Ferguson--Dirichlet process on $(\mathbb{X},\mathscr{X})$, with parameter-measure $\alpha$, if for every $k=1, 2, \dots$ and measurable partition $\{B_1, \dots, B_k\}$ of $\mathbb{X}$,
\begin{equation*}
	(\tilde{P}(B_1), \dots, \tilde{P}(B_k)) \sim \mathscr{D}(\alpha(B_1), \dots, \alpha(B_k)).
\end{equation*}
\end{definition}
The existence of the Ferguson--Dirichlet process follows from Lemma 1 in \cite{MR350949}, combined with the instrumental use of de~Finetti's representation theorem \citep{regazzini_petris92,regazzini1996}; see also Theorem 3.1 in \cite{ghosal2017}. The key ingredients of the proof are the additivity property \eqref{prop:additivity} and the observation that, for $k=1$, one has $B_1 = \mathbb{X}$ and $\tilde P(\mathbb{X}) = 1$ almost surely, by the definition of the Dirichlet distribution.

Note that Definition \ref{DefinitionDirichletProcess1} implies that $\tilde P$ is finitely additive. However, Proposition 2 in \cite{MR350949} shows that if $\alpha$ is $\sigma$--additive, then the Ferguson--Dirichlet process $\tilde P$ is $\sigma$-additive as well. Henceforth, we assume that $\alpha$ is $\sigma$-additive.

\subsection{Discreteness of the Ferguson--Dirichlet process}

In addition to its two fundamental properties, conjugacy and full support, the Ferguson--Dirichlet process prior also possesses another remarkable feature: it selects, almost surely, discrete distributions. This property was proved by \cite{blackwell1973}, implicitly derived  by \cite{MR350949} through the alternative construction presented in Section~\ref{Section_3.1}, and later proved using different techniques by \cite{berk1979}. The discreteness of the Ferguson--Dirichlet process plays a central role in many subsequent developments, in particular related to random partitions and clustering.

\begin{proposition}[\citealp{blackwell1973,MR350949}]
If $\tilde{P} \sim \mathscr{D}(\alpha)$, then $\tilde P$ is discrete, almost surely, that is  
\begin{equation}
	\label{eq:discrete_dp}
	\tilde{P} = \sum_{i \geq 1} \tilde{p}_i \delta_{Z_i},
\end{equation}
for some sequences $(\tilde p_i)_{i\ge 1}$ and $(Z_i)_{i\ge 1}$ of random probability weights and atoms, respectively. 
\end{proposition}
This result implies that, if $(X_i)_{i\ge 1}$ is an exchangeable sequence whose de~Finetti measure is a Ferguson--Dirichlet process prior $\mathscr{D}(\alpha)$, then $\mathbb{P}(X_i=X_j)>0$ for any $i\ne j$. Hence, for any $n\ge 1$, the first $n$ elements $X^{(n)}= (X_1, \dots, X_n)$ exhibit $K_n\le n$ distinct values $X_1^* , \dots, X_{K_n}^*$, with respective frequencies $n_1, \dots, n_{K_n}$. This induces a random partition of $[n] = \{1, \ldots, n\}$, where $i$ and $j$ are clustered together, denoted by $i\sim j$, if and only if $X_i=X_j$. The probability distribution of such a random partition is known as the Ewens sampling formula \citep{ewens1972}, which was also independently obtained by \cite{antoniak1974}, a student of T.S. Ferguson, who studied mixtures of Ferguson--Dirichlet processes. See also \cite{crane2016} and \cite{tavare2021}. The theory of random partitions, beyond the Ferguson--Dirichlet process case, has been central to many recent advances in Bayesian Nonparametrics. For a comprehensive account of the theory of exchangeable random partitions, see \cite{pitman2006}, and for their role in Bayesian nonparametric inference, see \cite{MR2730661} and \cite{de2013gibbs}. This line of research is closely connected to the literature on \emph{species sampling models} \citep{MR1481784}. The species metaphor applies to settings where a population is viewed as partitioned into infinitely many groups, each representing a distinct species. In this framework, $\tilde{p}_i$ in \eqref{eq:discrete_dp} represents the proportion of the population belonging to the species labeled by $Z_i$, while in a sample of size $n$, $K_n$ denotes the number of distinct observed species, with labels $X_j^*$ and frequencies $n_j$, for $j=1, \ldots, K_n$. For a review of species sampling models and their applications, particularly in predictive inference, see \cite{de2013gibbs}. More recently, \cite{flpr2025} proposed a multivariate version of species sampling models that generalizes several key properties of Pitman's formulation to the setting of non-exchangeable data. This framework also provides a natural foundation for studying distributional properties of models for partially exchangeable data, to be discussed in Section~\ref{Section5}.

The discreteness of realizations from a Ferguson--Dirichlet process makes it unsuitable for modeling continuous distributions directly. However, this property becomes highly valuable when defining priors on density functions through kernel mixtures of the form
\begin{equation*}
	\tilde{f}(x)=\int_{\mathbb{Y}}k(x,y)\,\tilde P(d y),
\end{equation*}
where $\tilde P \sim \mathscr{D}(\alpha)$. This approach, introduced in the pioneering works of \cite{lo1984} and \cite{ferguson1983}, has become one of the most successful applications of the Ferguson--Dirichlet process. Its enduring popularity and widespread use in density estimation and clustering are also due to \cite{escobar1995}, who were the first to develop a Gibbs sampling scheme  making Ferguson--Dirichlet process mixture models concretely implementable. See also \cite{mac98}, \cite{neal2000} and \cite{walker2007} for further key computational developments. In such mixture models, different choices of the prior on $\tilde P$ induce different random partitions of the latent variables and, consequently, different behavior for the number and relative sizes of the occupied mixture components.

\subsection{Alternative nonparametric priors}\label{sec:NIG}

It is natural to ask whether other distributions on the simplex can be used to define a nonparametric prior, thereby extending the construction in terms of finite-dimensional distributions beyond the Ferguson--Dirichlet process. A close inspection of Ferguson's construction reveals that the key lies in identifying distributions on the simplex that satisfy an additivity property analogous to that of the Dirichlet distribution \eqref{prop:additivity}. Once such a distribution is available, Ferguson's approach can be used as a blueprint for defining a corresponding alternative nonparametric prior.

Defining a distribution directly on the simplex is not straightforward. However, an entire class of distributions on the simplex can be obtained by mimicking the well-known construction of the Dirichlet distribution:
\smallskip

\noindent \textsc{Normalization.}
For any $i=1,\ldots, K$, let $\alpha_i>0$ and assume that
$Z_i \stackrel{\mbox{\scriptsize{\rm ind}}}{\sim} \text{\rm Ga}(\alpha_i, 1)$ are independent gamma random variables. If
\begin{equation}\label{eq:norm_gamma}
	Y_i := \frac{Z_i}{\sum_{j=1}^{K}Z_j} 
\end{equation}
for any $i=1,\ldots, K$, then
\begin{equation*}
	(Y_1, \dots, Y_K) \sim \mathscr{D}(\alpha_1, \dots, \alpha_K).
\end{equation*}
\smallskip

By replacing the gamma random variables in \eqref{eq:norm_gamma} with any other positive random variables, one can readily obtain a broad class of distributions on the simplex. For example, \cite{fhp2011} considered distributions on the simplex obtained by normalizing positive infinitely divisible random variables. However, the resulting distributions generally do not possess the additivity property.

A noteworthy special case is the normalized inverse-Gaussian distribution, introduced by \cite{MR2236441} and used by the authors to define the normalized inverse-Gaussian (N-IG) process. We briefly outline the key steps of its construction, in close analogy with Ferguson's approach.

Recall that a random variable $Z$ has the inverse-Gaussian distribution with shape parameter $\alpha>0$ and scale parameter $\gamma>0$, denoted by $Z \sim \text{IG}(\alpha, \gamma)$, if it admits the density
\begin{equation*}
	f_Z(z) = \frac{\alpha}{\sqrt{2\pi}} z^{-\frac{3}{2}} \:
	\mathrm{e}^{-\frac{1}{2} ( \frac{\alpha^2}{z} + \gamma^2 z ) + \gamma\alpha}
	\:\mathds{1}_{(0,+\infty)}(z).
\end{equation*}
A remarkable property of inverse-Gaussian random variables is their closure under convolution when they share the same scale parameter but may have different shape parameters. This property, also enjoyed by gamma random variables, can be stated as follows: if $Z_i\stackrel{\mbox{\scriptsize ind}}{\sim} \mbox{IG}(\alpha_i,\gamma)$ for $i=1,\ldots,K$, then
\begin{equation}\label{IGdistributionAdditivity}
	\sum_{i=1}^K Z_i\sim\mbox{IG}(|\bm{\alpha}|,\,\gamma)
\end{equation}
with $\bm{\alpha}:=(\alpha_1, \dots, \alpha_K)$ and $|\bm{\alpha}|=\sum_{i=1}^K \alpha_i$. 	

\begin{figure}[t]
	\centering
	\includegraphics[width=0.5\textwidth]{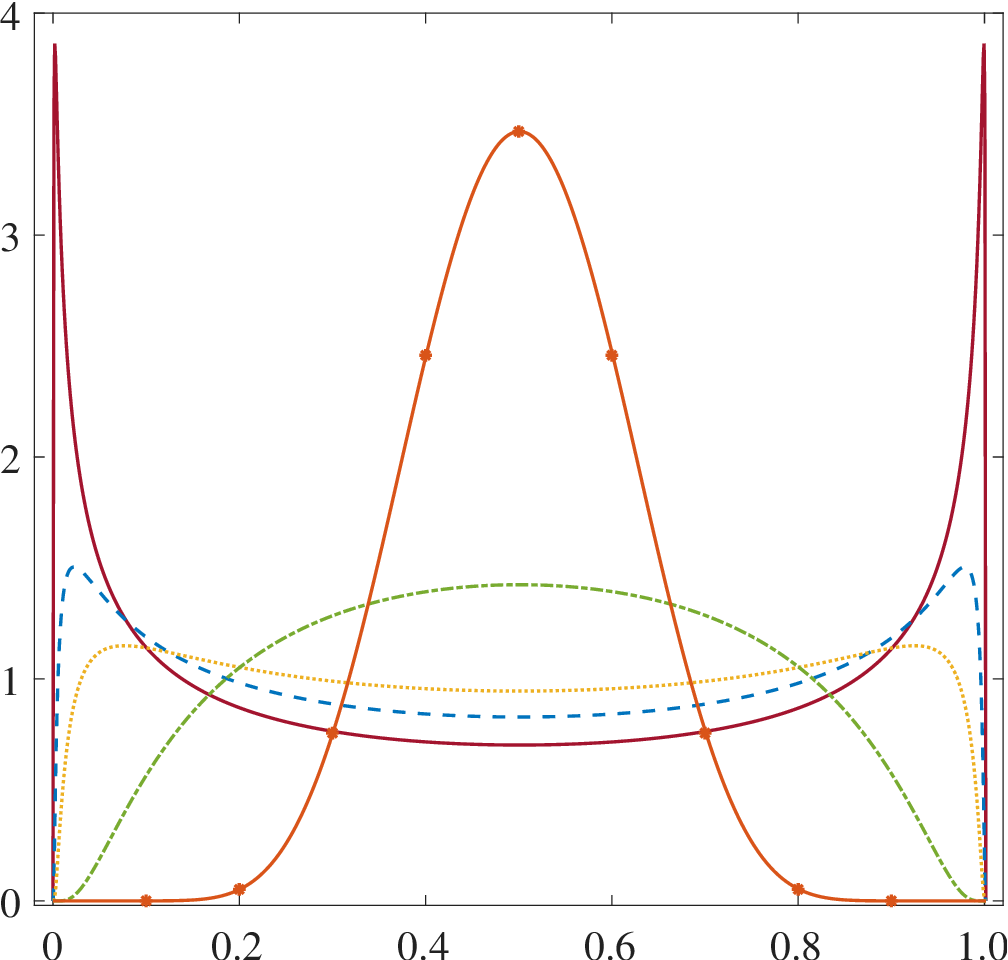}
	\caption{Normalized inverse Gaussian densities on $(0,1)$ with parameters $(\alpha_i, \alpha_i)$. \textcolor{mycolor1}{NIG$(.0583, .0583)$\rule[.5ex]{1.7em}{.8pt}}; \textcolor{mycolor2}{NIG$(.198, .198)$\rule[.5ex]{0.5em}{.8pt}\,\,\rule[.5ex]{0.5em}{.8pt}\,\,\rule[.5ex]{0.5em}{.8pt}};\textcolor{mycolor3}{NIG$(.354, .354)$\rule[.5ex]{0.1em}{.8pt}\,\rule[.5ex]{0.1em}{.8pt}\,\rule[.5ex]{0.1em}{.8pt}\,\rule[.5ex]{0.1em}{.8pt}\,\rule[.5ex]{0.1em}{.8pt}\,\rule[.5ex]{0.1em}{.8pt}}; \textcolor{mycolor4}{NIG$(1.23, 1.23)$\rule[.5ex]{0.5em}{.8pt}\,\rule[.5ex]{0.1em}{.8pt}\,\rule[.5ex]{0.5em}{.8pt}\,\rule[.5ex]{0.1em}{.8pt}\,\rule[.5ex]{0.5em}{.8pt}}; \textcolor{mycolor5}{NIG$(9.07, 9.07)$\rule[.5ex]{0.8em}{.8pt}$*$\rule[.5ex]{0.8em}{.8pt}}.}
	\label{FigureNIGDensityK2}
\end{figure}

\begin{figure}
	\centering
	\begin{subfigure}[t]{0.4\textwidth}
		\includegraphics[width=1\textwidth]{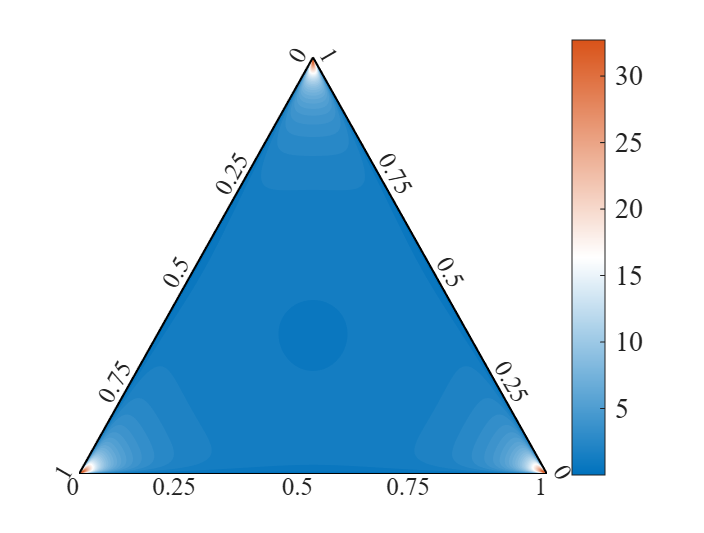}
		\caption{N-IG$(.183,.183,.183)$.}
	\end{subfigure}
	\enskip
	\begin{subfigure}[t]{0.4\textwidth}
		\includegraphics[width=1\textwidth]{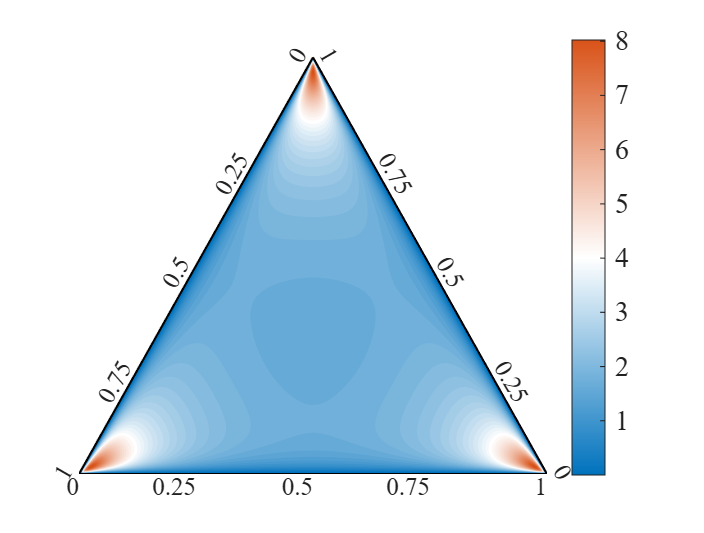}
		\caption{N-IG$(.346,.346,.346)$.}
	\end{subfigure}
	
	\begin{subfigure}[t]{0.4\textwidth}
		\includegraphics[width=1\textwidth]{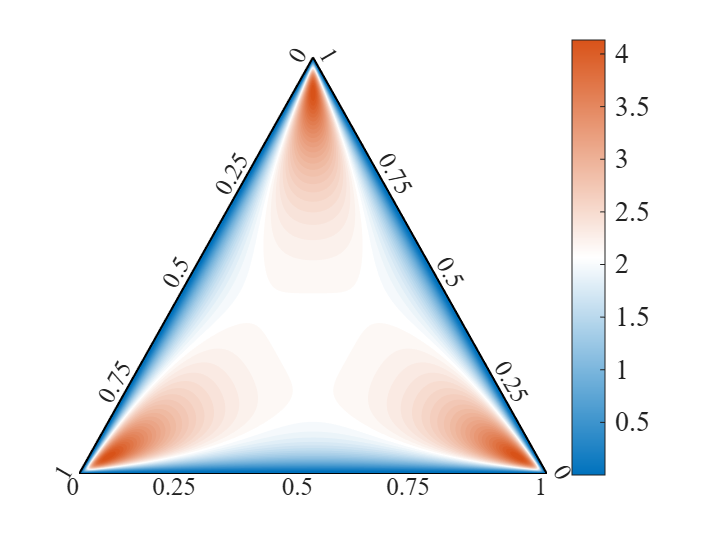}
		\caption{N-IG$(.52,.52,.52)$.}
	\end{subfigure}
	\enskip
	\begin{subfigure}[t]{0.4\textwidth}
		\includegraphics[width=1\textwidth]{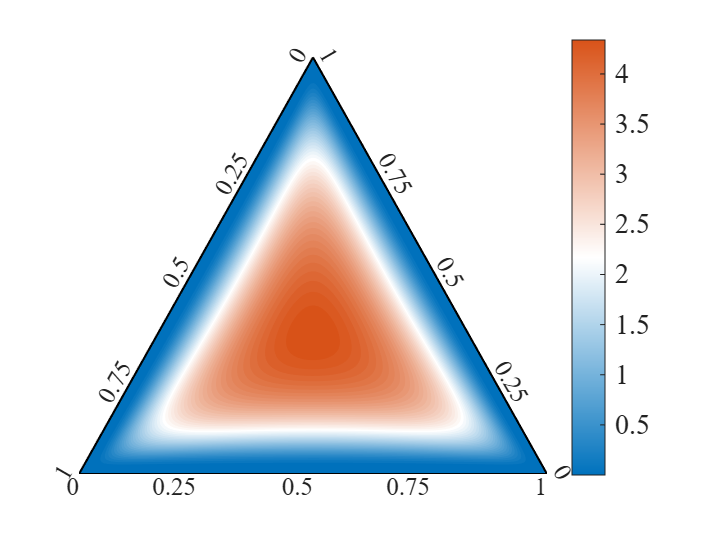}
		\caption{N-IG$(1.45,1.45,1.45)$.}
	\end{subfigure}
	\caption{Heatmaps for the normalized inverse Gaussian density on $\Delta_2$ corresponding to various parameter settings.}
	\label{FigureNIGDensityK3}
\end{figure}

The normalized inverse-Gaussian distribution is then obtained via normalization like the Dirichlet distribution in \eqref{eq:norm_gamma}: Let $Z_i \stackrel{\mbox{\scriptsize ind}}{\sim} \text{IG}(\alpha_i, 1)$, with $\alpha_i > 0$ for $i=1,\ldots,K$, and define the random vector 
\begin{equation}\label{NIGNormalizationConstruction}
	(Y_1, \dots, Y_K) := \left( \frac{Z_1}{\sum_{i=1}^{K}Z_i} ,\dots, \frac{Z_K}{\sum_{i=1}^{K}Z_i} \right)
\end{equation}
Then $(Y_1,\ldots,Y_K)$ has the normalized inverse Gaussian distribution with parameters $\alpha_1,\ldots,\alpha_K$, denoted by
\begin{equation*}
	(Y_1, \dots, Y_K) \sim \text{N-IG}(\alpha_1, \dots, \alpha_K).
\end{equation*}
The density of $(Y_1, \dots, Y_{K-1})$ equals
\begin{equation*}
	\begin{split}
		f_K(\bm{y}; \bm{\alpha})
		&= 
		\frac{\mathrm{e}^{|\bm{\alpha}|}\prod_{i=1}^{K}\alpha_i}
		{2^{\frac{K}{2}-1} \pi^{\frac{K}{2}}}
		\:\mathcal{K}_{-K/2}\Big(\sqrt{
			\mathcal{A}_K(\bm{y})}\Big)  y_1^{-\frac{3}{2}} \dots y_{K-1}^{-\frac{3}{2}} \left( 1 - |\bm{y}|  \right)^{-\frac{3}{2}} [\mathcal{A}_K(\bm{y})]^{-\frac{K}{4}} \:\mathds{1}_{\Delta_{K-1}}(\bm{y}),
	\end{split}
\end{equation*}
with $\mathcal{K}_{-\frac{K}{2}}(\cdot)$ being the modified Bessel function of the third type, and
\begin{equation*}
	\mathcal{A}_K(\bm{y}; \bm{\alpha}) := \sum_{i=1}^{K-1} \frac{\alpha_i^2}{y_i} + \frac{\alpha_K^2}{1 -|\bm{y}|}. 
\end{equation*}
Figures~\ref{FigureNIGDensityK2} and~\ref{FigureNIGDensityK3} illustrate the normalized inverse Gaussian distribution on $(0,1)$ and $\Delta_2$, respectively, under different parameter settings. These plots display interesting behaviors near the boundaries, especially for small values of the $\alpha_i$'s,  quite different from the Dirichlet distribution. To facilitate the comparison, we choose the parameters of the N-IG distributions so that their means and variances coincide with those of the Dirichlet distributions shown in Figure~\ref{FigureDirichletDensityK2} and Figure~\ref{FigureDirichletDensityK3}. Recall that if $(Y_1, \dots, Y_K) \sim \mathscr{D}(\alpha_1, \dots, \alpha_K)$, then $\mathbb{E}(Y_k) = \alpha_k / |\bm{\alpha}|:=p_k^0$ and $\text{Var}(Y_k) =(1 + |\bm{\alpha}|)^{-1} p_k^0 (1- p_k^0)$, for any $k = 1, \ldots , K$. On the other hand, for $(Y_1, \dots, Y_K) \sim \mathrm{N\mbox{-}IG}(\bar{\alpha}_1, \dots, \bar{\alpha}_K)$, $\mathbb{E}(Y_k) = \bar{\alpha}_k / |\bar{\bm{\alpha}}|:=\bar p_k^0$ and $\text{Var}(Y_k) =  |\bar{\bm{\alpha}}|^2 \mathrm{e}^{|\bar{\bm{\alpha}}|} \Gamma(-2, |\bar{\bm{\alpha}}|) \bar p_k^0 (1- \bar p_k^0)$, where $\Gamma(x,a)$ denotes the incomplete gamma function. See \cite{MR2236441} for details. Hence the means coincide when  $p_k^0 = \bar p_k^0$, whereas, for a given $|\bm{\alpha}|$, matching the variances amounts to solving $(1 + |\bm{\alpha}|)^{-1} = |\bar{\bm{\alpha}}|^2 \mathrm{e}^{|\bar{\bm{\alpha}}|} \Gamma(-2, |\bar{\bm{\alpha}}|)$ with respect to $|\bar{\bm{\alpha}}|$.

By virtue of \eqref{IGdistributionAdditivity}, one can show the additivity of the N-IG distribution: 

\smallskip
\noindent \textsc{Additivity.}
If $(Y_1,\ldots,Y_K)\sim\mbox{\rm N-IG}(\alpha_1,\ldots,\alpha_K)$, and $r_1,\ldots,r_{l-1}$ are positive integers such that $r_1<r_2<\,\cdots\,<r_{l-1}<K$, then
\begin{equation}\label{prop:additivity_nig} 
\left( \sum_{i=1}^{r_1}Y_i , \sum_{i=r_1+1}^{r_2}Y_i , \dots , \sum_{i=r_{l-1}+1}^{K}Y_i \right)
\sim \textstyle \mbox{\rm N-IG} \left( \sum_{i=1}^{r_1}\alpha_i , \sum_{i=r_1+1}^{r_2}\alpha_i , \dots , \sum_{i=r_{l-1}+1}^{K}\alpha_i \right) .
\end{equation}

With an additive distribution on the simplex \eqref{prop:additivity_nig} in hand, the essential ingredient is now in place. Following Ferguson's foundational construction, the existence of a random probability measure with N-IG finite-dimensional distributions can be established; this is the normalized inverse-Gaussian process.

\begin{definition}[N-IG process, \citealt{MR2236441}]
Let $\alpha$ be a non-null finite  measure (non-negative and finitely additive) on $(\mathbb{X},\mathscr{X})$. We say $\tilde P$ is a normalized inverse Gaussian (N-IG) process on $(\mathbb{X},\mathscr{X})$ with parameter-measure $\alpha$ if for every $k\ge 1$ and measurable partition $\{B_1, \dots, B_k\}$ of $\mathbb{X}$, one has
\begin{equation*}
	(\tilde P(B_1), \dots, \tilde P(B_k)) \sim \text{N-IG}(\alpha(B_1), \dots, \alpha(B_k)).
\end{equation*}
\end{definition}

A natural question is whether the N-IG process also satisfies the two desiderata, large support and analytical tractability. Regarding the former, it can be shown that, with respect to the weak topology, its support is the set of all probability measures whose support is contained in that of the parameter measure $\alpha$. Regarding tractability, although the N-IG process is not conjugate \citep{MR2255112}, it remains analytically tractable: predictive distributions \citep{MR2236441} and a posterior characterization \citep{MR2508332} are available. The N-IG process is therefore a viable alternative nonparametric prior.

As noted previously, distributions on the simplex obtained by normalizing positive random variables generally do not yield additivity. Besides the gamma and inverse Gaussian cases (the latter including the notable $1/2$-stable special case), the only other known example that does is a generalization of the gamma law used in \cite{lmp05} to define a generalized Ferguson--Dirichlet process, which exhibits interesting connections to Lauricella hypergeometric functions \citep{lauricella1893}. Outside normalization schemes for constructing distributions on the simplex, only one additive family is currently known. Determined by \cite{carlton2002}, another student of T.S.~Ferguson, it coincides with the finite-dimensional distributions of the Pitman-Yor process \citep{MR1434129} with parameters $\theta$ and $\sigma = 1/2$. Identifying a broad class of additive distributions on the simplex, or at least further new special cases, remains an open problem.

\section{Construction via the normalization of completely random measures}\label{Section3}

In this section, we review, from a modern perspective, the second construction of the Ferguson--Dirichlet process, introduced in \cite{MR350949}, as the normalization of a suitable gamma process. Like the construction based on finite-dimensional distributions, this approach provides a natural gateway to generalizations obtained by replacing the gamma process with other processes having positive independent increments. Moreover, unlike the finite-dimensional approach, it allows the definition of an entire class of nonparametric priors through the normalization of general completely random measures (CRMs), culminating in the formulation of \emph{normalized random measures with independent increments} \citep{MR1983542}.

\subsection{The normalized gamma process}\label{Section_3.1}

We describe the Ferguson--Dirichlet process as the normalization of a gamma process, following the intuition outlined by Ferguson and showing how this yields an equivalent and insightful formulation.
Definition~\ref{DefinitionDirichletProcess1} links the Ferguson--Dirichlet process to the (finite-dimensional) Dirichlet distribution, while the latter arises from the normalization of independent gamma random variables as in~\eqref{eq:norm_gamma}. Combining these two facts, if $Z_1, \ldots, Z_K$ are independent random variables with $Z_i \sim \mathrm{Ga}(\alpha(B_i), 1)$, then
\begin{equation*}
	(\tilde P(B_1), \ldots, \tilde P(B_K))
	\stackrel{d}{=}
	\left(\frac{Z_1}{\sum_{i=1}^{K} Z_i}, \ldots, \frac{Z_K}{\sum_{i=1}^{K} Z_i}\right).
\end{equation*}
This observation suggests that the Ferguson--Dirichlet process can be defined in an analogous way by replacing the finite collection of independent gamma variables with a \emph{gamma process}, a process with independent gamma-distributed increments, and normalizing it by its total mass. For the normalization to be well defined, the total mass of the gamma process must be finite almost surely. As noted by \cite{MR350949}:

\begin{quote}
	\textsl{The basic idea is that since the Dirichlet distribution is definable as the joint distribution of a set of independent gamma variables divided by the sum, so also should the Dirichlet process be definable as a gamma process with independent increments divided by the sum.}
\end{quote}
This insight leads to an equivalent definition of the Ferguson--Dirichlet process as the normalization of a process with independent increments. 

Following the notation in \cite{MR350949}, consider a gamma process 
\begin{equation*}
	Z_t \in \mathscr{G}(\alpha((0,t]),1), \quad t \geq 0 ,
\end{equation*}
with  finite non-null parameter measure $\alpha$, that is, a process with independent increments and almost surely non-decreasing paths such that (a) $\mathbb{P}(Z_0=0)=1$ and (b) $Z_t - Z_s \sim \mathrm{Ga}(\alpha((s,t]),1)$ for any $0 \le s < t < \infty$. Note that $Z_t$ does not have stationary increments, and hence is not a gamma subordinator \citep{bertoin}. Indeed, a gamma subordinator would have stationary increments, which would require $\alpha$ to be proportional to the Lebesgue measure and would imply $\lim_{t\to\infty}Z_t=\infty$. Set $F_0(t) = \alpha((0,t]) / \alpha(\mathbb{R}^+)$. In \cite{MR373022} it is shown that this gamma process admits the following explicit series representation
\begin{equation}\label{eq:gamma_series}
	Z_t = \sum_{i=1}^\infty J_i\,\mathds{1}_{[0,F_0(t))}(U_i), \quad J_1>J_2>\ldots>0 ,
\end{equation}
where $U_i\stackrel{\mathrm{iid}}{\sim}\mathrm{Unif}(0,1)$, the positive jumps $J_i\in\mathbb R^+$ are arranged in decreasing order, and, for any $x>0$,
\begin{equation*}
	\mathbb{P}(J_1 \le x)
	= \exp\!\left(-\alpha(\mathbb{R}^+)\int_x^{\infty} \mathrm{e}^{-s}s^{-1}\,ds\right),
\end{equation*}
while, for any $j \ge 2$ and $0<x_j<x_{j-1}$,
\begin{equation*}
\mathbb{P}[ \, J_j\le x_j\, \mid \, J_1=x_1,\ldots,J_{j-1}=x_{j-1}]
=
\exp\Big(-\alpha(\mathbb{R}^+)\,\int_{x_{j}}^{x_{j-1}}\mathrm{e}^{-s}\, s^{-1}\:ds \Big).
\end{equation*}
The sequences $(U_i)_{i \ge 1}$ and $(J_i)_{i \ge 1}$ are independent. See also Theorem~1 in Section 4 of \cite{MR350949}. This elegant representation highlights the discreteness of the gamma process and provides an intuitive constructive view of its sample paths.

Building on this representation, Theorem~2 in Section 4 of \cite{MR350949} provides an equivalent definition of the Ferguson--Dirichlet process expressed in terms of the random distribution function $t \mapsto \tilde{F}(t) = \tilde{P}((0,t])$. 

\begin{definition}[Dirichlet process, \citealt{MR350949}]\label{DefinitionDirichletProcess2}
	Let $Z_t \in \mathscr{G}(\alpha((0,t]),1)$ be a gamma process, and let $Z_{\infty} = \lim_{t \to \infty} Z_t \in \mathscr{G}(\alpha(\mathbb{R}^+),1)$. 
	Then $\tilde P \sim \mathscr{D}(\alpha)$ (or equivalently $\tilde F \sim \mathscr{D}(\alpha)$) if
	\begin{equation*}
		\tilde{F}(t) = \frac{Z_t}{Z_{\infty}}, \quad t \ge 0.
	\end{equation*}
\end{definition}

This alternative definition yields an interesting representation of the Ferguson--Dirichlet process as an almost surely discrete random probability measure. Let $P_0:=\alpha/\alpha(\mathbb{R}^+)$, so that $F_0(t)=P_0((0,t])$. Consider a sequence of i.i.d. random variables $(V_i)_{i\ge 1}$, independent of $(J_i)_{i\ge 1}$, with common distribution $P_0$. Defining $\tilde p_i=J_i/Z_\infty$, for $i\ge 1$, we obtain
\begin{equation*}
	\tilde{P}=\sum_{i\ge 1} \tilde p_i\: \delta_{V_i},
\end{equation*}
which is a Ferguson--Dirichlet process with parameter measure $\alpha$.

\subsection{Completely random measures}

Ferguson's normalization approach prompts a broader question: can the gamma process be replaced by a more general class of processes, while retaining analytic tractability of the resulting nonparametric prior? A.F.M.~Smith posed this challenge in his discussion of \cite{MR368264}:
 \begin{quote}
	\textsl{So far as I know, Ferguson's approach using the Dirichlet process is the only procedure available which enables one to work through, analytically, the usual prior to posterior steps. Appendix 2 points out that, in a sense, Ferguson's Dirichlet process is a special case of a rather more general class of processes. The question of interest to a Bayesian statistician is whether there are any other processes in this class which are tractable. The same kind of question arises in connection with the generation of distributions by means of subordinators.}
\end{quote}
The answer is affirmative. The ``more general class'' is known as \emph{normalized random measures with independent increments} (NRMIs) \citep{MR1983542} in the modern Bayesian nonparametric literature. The idea is to replace the gamma process in Definition~\ref{DefinitionDirichletProcess2} with a general completely random measure (CRM) and normalize it; the Ferguson--Dirichlet process is recovered when the CRM is gamma.

Let $\mathbb{M}$, equipped with its Borel $\sigma$--algebra $\mathscr{M}$,  denote the space of boundedly finite measures on $\mathbb{X}$, that is, $\mu(E) < \infty$ for every bounded $E \in \mathscr{X}$; see \cite{daley08} for details.
\begin{definition}[Completely random measure, \citealt{MR210185, MR1207584}]
	Let $\tilde \mu$ be a random element defined on some probability space $(\Omega, \mathscr{F}, \mathbb{P})$ and taking values in $(\mathbb{M},\mathscr{M})$. If for any $d \geq 2$ and pairwise disjoint sets $A_1, \ldots, A_d$ in $\mathscr{X}$, the random variables $\tilde{\mu}(A_1), \dots , \tilde{\mu}(A_d)$ are mutually independent, then $\tilde\mu$ is said to be a \emph{completely random measure} (CRM) on $\mathbb{X}$.
\end{definition}

Kingman's decomposition theorem \citep{MR210185} represents any completely random measure (CRM) as the sum of a deterministic drift, a random component with jumps at fixed locations, and a random component with jumps at random locations. For our purposes, it is enough to consider CRMs without drift and without fixed atoms. In this case the CRM is almost surely discrete and admits the representation 
\begin{equation}\label{CRMdef}
	\tilde{\mu} \;=\; \sum_{i=1}^{\infty}\tilde{J}_i \delta_{V_i},
\end{equation}
where $(\tilde J_i)_{i\ge 1}$ are positive random jumps and $(V_i)_{i\ge 1}$ are $\mathbb{X}$-valued random locations. The pairs $\{(\tilde J_i,V_i): i\ge 1\}$ are the points of a Poisson random measure with mean measure $\nu$ on $\mathbb{R}^+\times\mathbb{X}$ such that
\begin{equation}
	\label{eq:intensity_condition}
	\int_B\int_{\mathbb{R}^+}\min\{v,1\}\,\nu(dv,dx)<+\infty
\end{equation} 
for any bounded $B$ in $\mathscr{X}$. Under \eqref{eq:intensity_condition}, the Laplace functional of $\tilde{\mu}$ admits the following representation \citep{MR1207584}.
\begin{proposition}[Laplace functional of a CRM] \label{prop:laplace}
	A CRM $\tilde{\mu}$ is uniquely characterised by the Laplace functional transform
	\begin{equation*}
		\mathbb{E}\left( \exp\left( -\int_{\mathbb{X}} g(x)\tilde{\mu}(dx) \right) \right) = \exp\left( -\int_{\mathbb{R}^+ \times \mathbb{X}} (1 - \mathrm{e}^{-vg(x)}) \nu(dv,dx) \right) ,
	\end{equation*}
	for any bounded $g:\mathbb{X}\to\mathbb{R}^+$. The measure $\nu$ satisfying the integrability condition \eqref{eq:intensity_condition} is the L\'evy intensity of $\tilde\mu$.
\end{proposition}
\noindent Hence, any measure satisfying \eqref{eq:intensity_condition} identifies a CRM on $\mathbb{X}$. Borrowing the terminology of Ferguson, we distinguish two cases:
\begin{itemize}
	\item If $\nu(dv, dx) = \rho(dv)\alpha(dx)$, the jump sizes and locations are independent, and $\tilde{\mu}$ is a \emph{homogeneous CRM}.
	\item If $\nu(dv, dx) = \rho_x(dv)\alpha(dx)$, the distribution of the jump size depends on the location, and $\tilde{\mu}$ is termed a \emph{non-homogeneous CRM}.
\end{itemize}
Relevant examples widely used in Bayesian Nonparametrics, identified by their L\'evy intensities, include:
\begin{enumerate}
	\item Gamma CRM \citep{MR350949}: for any $\tau>0$
	\begin{equation*}
		\nu(dv,dx) = \frac{\mathrm{e}^{-\tau v}}{v}dv \alpha(dx) ;
	\end{equation*}
	\item $\sigma$-stable CRM \citep{MR368264}: for any $\sigma\in(0,1)$
	\begin{equation*}
		\nu(dv,dx) = \frac{\sigma}{\Gamma(1-\sigma) v^{1+\sigma}}dv \alpha(dx) ;
	\end{equation*}
	\item Inverse--Gaussian CRM \citep{MR2236441}: for any $\tau>0$ 
	\begin{equation*}
		\nu(dv,dx) = \frac{\mathrm{e}^{-\tau v}}{\sqrt{2\pi}v^{3/2}}dv \alpha(dx) ;
	\end{equation*}
	\item Generalized gamma \citep{MR1747450}, also known as tempered stable \citep{rosinski}, CRM: for any $\tau>0$ and $\sigma\in(0,1)$
	\begin{equation*}
		\nu(dv,dx) = \frac{\sigma \mathrm{e}^{-\tau v}}{\Gamma(1 - \sigma) v^{1+\sigma}}dv \alpha(dx) ;
	\end{equation*}
	\item Extended gamma CRM \citep{MR606619}: if $\beta:\mathbb{X}\to\mathbb{R}^+$ is such that $\int_{B}\log(1+\beta(x))\alpha(dx)<+\infty$ for any bounded $B$ in $\mathscr{X}$
	\begin{equation*}
		\nu(dv,dx) = \frac{\mathrm{e}^{-\beta(x) v}}{v}dv \alpha(dx) ;
	\end{equation*}
	\item Simple homogeneous CRM \citep{ferguson1979}: for any $\tau>0$
	\begin{equation*}
		\nu(dv,dx) = \frac{1 - \mathrm{e}^{-\tau v}}{1 - \mathrm{e}^{-v}} \frac{\mathrm{e}^{-v}}{v} dv \alpha(dx) .
	\end{equation*}
\end{enumerate}

Because of the infinite series representation \eqref{CRMdef}, it is impossible to sample realizations of a CRM exactly, and one must rely on approximations. Various methods have been studied in the literature, and the most efficient approach is based on the \emph{Ferguson and Klass representation}, which generalizes the gamma process series in \eqref{eq:gamma_series}; see Section~4 of \cite{MR350949} and \cite{MR373022}. An important practical advantage of this representation is that it generates the jumps in decreasing order \citep{MR1782492}.

We briefly outline the extension of the gamma process representation to any homogeneous CRM. For a homogeneous CRM with $\nu(dv,dx)=\rho(dv)\alpha(dx)$, let $J_1>J_2>\cdots$ denote the decreasing rearrangement of the jumps $(\tilde J_i)_{i\ge1}$ in \eqref{CRMdef}. These ranked jumps can be sampled by $J_i=\rho^{\leftarrow}(\Gamma_i)$, where $\rho^{\leftarrow}(w):=\inf\{ x>0 \mid \: \alpha(\mathbb{X}) \int_x^\infty \rho(dv)\leq w \}$, $\Gamma_k:=\sum_{j=1}^{k}E_j$, and $E_j\sim\text{Exp}(1)$ are i.i.d. exponential random variables with mean $1$. This method is known as the inverse-L\'evy measure method \citep{MR1833707}. In principle, this method applies to any CRM for which the inverse-L\'evy measure $\rho^{\leftarrow}$ can be evaluated. In many Bayesian nonparametric applications it is also necessary to sample the total mass of the remaining jumps $T_{K+1}\;=\;\sum_{i=K+1}^{\infty}J_i$. For instance, in finite approximations to mixture models, the likelihood depends on the $K$ largest jumps $J_1,\ldots,J_K$ and on the rest only through $T_{K+1}$. Sampling $T_{K+1}$ amounts to simulating a truncated subordinator, for which no general algorithm is available, although effective solutions exist for several widely used priors; see \cite{MR4530320} and \cite{zhang2024posterior}.

\subsection{Normalized CRMs}\label{sec:NRMI}

The construction of the Ferguson--Dirichlet process in  Definition~\ref{DefinitionDirichletProcess2} hinges on $Z_\infty$ being almost surely finite and strictly positive. To extend this approach beyond the Ferguson--Dirichlet process, one must ensure that the following two conditions are satisfied:
\begin{itemize}
	\item[(H1)] \emph{Finiteness:} $\tilde{\mu}(\mathbb{X})<\infty$ almost surely; for a homogeneous CRM, this amounts to $\alpha(\mathbb{X})<\infty$.
	\item[(H2)] \emph{Infinite activity:} $\tilde{\mu}(\mathbb{X})>0$ almost surely; for a homogeneous CRM with $\alpha(\mathbb X)>0$, this is equivalent to $\rho(\mathbb{R}^+) = \infty$.
\end{itemize}

We now introduce normalized random measures with independent increments, a general framework for constructing random probability measures. On the one hand, it allows one to recover the Ferguson--Dirichlet process as the normalization of a gamma CRM, extending Ferguson's original setting $\mathbb{X}=\mathbb{R}^+$. On the other hand, it yields a broad class of nonparametric priors that generalize the Ferguson--Dirichlet process. 
\begin{definition}[NRMI, \citealp{MR1983542}]\label{DefinitionNRMI}
	Let $\tilde{\mu}$ be a CRM on $\mathbb{X}$ satisfying (H1) and (H2). Then the random probability measure
	\begin{equation*}
		\tilde{P} := \frac{\tilde{\mu}}{\tilde{\mu}(\mathbb{X})}
	\end{equation*}
	is called a \emph{normalized random measure with independent increments} (NRMI). 
\end{definition}

The term ``increments'' in Definition~\ref{DefinitionNRMI} reflects that the original formulation in \cite{MR1983542} considered the case $\mathbb{X}=\mathbb{R}$. Note that an NRMI is determined by the L\'evy intensity $\nu$ of its underlying CRM. For example, the Ferguson--Dirichlet process is recovered as a special case in which $\tilde{\mu}$ is a gamma CRM with finite base measure $\alpha$. 

NRMIs, similarly to the Ferguson--Dirichlet process, satisfy Ferguson's two \emph{desiderata}. Indeed, the support of an NRMI is the set of all probability measures whose support is contained  in the support of its parameter--measure $\alpha$. In particular, if $\operatorname{supp}(\alpha)=\mathbb{X}$, then the support coincides with the entire space $\mathscr{P}$ of probability measures on $\mathbb{X}$. Analytical tractability is attained both \emph{a priori} and \emph{a posteriori}. As to the former, prior moments of NRMIs, which are helpful for guiding prior specification, can be determined. The following result records second order moments, though higher order moments can also be obtained.

For simplicity, we henceforth focus on homogeneous NRMIs and set $\theta := \alpha(\mathbb{X})$ and $P_0 := \alpha/\alpha(\mathbb{X})$, so that $\alpha = \theta P_0$.

\begin{proposition}[\citealp{MR2255112}]\label{TheoremAprioriMoments}\label{prop:NRMI_mom}
	Let $\tilde{P}$ be a homogeneous NRMI with L\'evy intensity $\nu(dv,dx) = \rho(dv)\theta P_0(dx)$. Then for any $A, B \in \mathscr{X}$, 
	\begin{equation*}
		\begin{split}
			\mathbb{E}(\tilde{P}(A)) &= P_0(A) ,\\
			\textnormal{Var}(\tilde{P}(A)) &= P_0(A)P_0(A^c) \mathscr{I}_{\theta} ,
		\end{split}
	\end{equation*}
	and
	\begin{equation*}
		\textnormal{Cov}(\tilde{P}(A), \tilde{P}(B)) = \{P_0(A\cap B) - P_0(A)P_0(B)\} \mathscr{I}_{\theta}, \vspace{3pt}
	\end{equation*}
	with $\tau_i(u) = \int_{0}^{\infty} v^i \mathrm{e}^{-uv} \rho(dv)$ for $u>0$ and $i=1,2,\ldots$ and
	\begin{equation*}
		\mathscr{I}_{\theta} = \theta \int_{0}^{\infty} u \tau_2(u) \mathrm{e}^{-\theta\int (1 - \mathrm{e}^{-uv})\rho(dv)} du .
	\end{equation*}
\end{proposition}
Note first that the prior mean $\mathbb{E}\{\tilde P(A)\}$ coincides with $P_0(A)$ regardless of the choice of the L\'evy intensity $\rho$. Thus all homogeneous NRMIs with the same $P_0$ share the same first moment. Moreover, for homogeneous NRMIs the second-order structure is universal up to a multiplicative constant, which depends on the L\'evy intensity $\rho$, and thus on the specific NRMI. The \textit{a priori} moments of the Ferguson--Dirichlet process follow immediately as a special case of Proposition~\ref{TheoremAprioriMoments}.

\begin{proposition}[\citealp{MR350949}, Proposition 1] \label{prop:DP_mom}
	Let $\tilde{P} \sim \mathscr{D}(\alpha)$. Then for any $A, B \in \mathscr{X}$, 
	\begin{equation*}
		\begin{split}
			\mathbb{E}(\tilde{P}(A)) &= \frac{\alpha(A)}{\alpha(\mathbb{X})} = P_0(A) ,\\
			\textnormal{Var}(\tilde{P}(A)) &= \frac{P_0(A)P_0(A^c)}{\theta+1} ,\\
			\textnormal{Cov}(\tilde{P}(A), \tilde{P}(B)) &= \frac{P_0(A\cap B) - P_0(A)P_0(B)}{\theta+1} .
		\end{split}
	\end{equation*}
\end{proposition}
Some comments on the proofs are in order. Proposition~\ref{prop:DP_mom} was proved by Ferguson relying on the fact that the marginal distributions of a Ferguson--Dirichlet process are explicitly known \eqref{MarginalDistributionDirichletProcess}. Within the class of homogeneous NRMIs, explicit marginals are the exception rather than the rule. By contrast, the L\'evy intensity and the associated Laplace functional given in Proposition~\ref{prop:laplace} are available, and these are the tools used to establish Proposition~\ref{prop:NRMI_mom}. In the Appendix  we provide an outline of the proof to illustrate a key technique to study CRM--based Bayesian nonparametric models.

Analytical tractability of NRMIs is also preserved \emph{a posteriori}. At first sight this claim may seem at odds with \cite{MR2255112}, which shows that within homogeneous NRMIs conjugacy is unique to the Ferguson--Dirichlet process.
 
\begin{theorem}[\citealp{MR2255112}]\label{thm:conjugacy}
Let $\mathcal{Q}$ be the set of all probability distributions of homogeneous NRMIs. The posterior distribution of a homogeneous NRMI $\tilde P$, given the observations drawn as in \eqref{BayesianHierarhicalModel}, belongs to $\mathcal{Q}$ if and only if $\tilde P$ is a Ferguson--Dirichlet process.
\end{theorem}

Although Theorem~\ref{thm:conjugacy} rules out both parametric and structural conjugacy for homogeneous NRMIs other than the Ferguson--Dirichlet process, their posteriors remain analytically tractable. In fact, \cite{MR2508332} provide an explicit posterior characterization of  NRMIs. It is shown that there exists a latent random variable conditional on which NRMIs are conjugate, which can be seen as a \emph{conditional structural conjugacy}. This characterization yields workable posterior and predictive distributions and underpins efficient Gibbs and marginal samplers for NRMI mixture models; see, e.g., \cite{statsci_2013}.

Let $X^{(n)}:=(X_1, \dots, X_n)$ be a sample of size $n$ from an NRMI $\tilde{P}$, according to \eqref{BayesianHierarhicalModel}. Assume that there are $K_n$ distinct values in $X^{(n)}$, namely $X_1^*, \dots, X_{K_n}^*$, with frequencies $n_1, \dots, n_{K_n}$, respectively. Let $U_n$ be a positive random variable whose density, conditional on $X^{(n)}$, is
\begin{equation*}
	f(u \mid X^{(n)}) \propto u^{n-1} \mathrm{e}^{-\psi(u)} \prod_{i=1}^{K_n} \tau_{n_i}(u) ,
\end{equation*}
with
\begin{equation*}
	\psi(u) := \theta \int_{\mathbb{R}^+ } (1 - \mathrm{e}^{-uv}) \rho(dv) ,
\end{equation*}
and, for any $m\geq 1$,
\begin{equation*}
	\tau_m(u) := \int_{\mathbb{R}^+} s^m \mathrm{e}^{-us} \rho(ds) .
\end{equation*}
Conditionally on the data $X^{(n)}$ and the latent variable $U_n$, the posterior distribution $\tilde{P} \mid (X^{(n)}, U_n)$ is still in $\mathcal{Q}$. A key step to establish this result is the posterior characterization of the CRM $\tilde\mu$:

\begin{theorem}[\citealp{MR2508332}]
	\label{TheoremNRMIposterior}
	The posterior distribution of $\tilde{\mu}$, given $X^{(n)}$, is a mixture with respect to $U_n$, namely
	\begin{equation}
		\label{eq:posterior_general}
		\left( \tilde{\mu} \mid X^{(n)} , U_n \right) = \tilde{\mu}_{_{U_n}} + \sum_{i=1}^{K_n}J_i^{(U_n)} \delta_{X_i^*} ,
	\end{equation}
	where, for any $u>0$,
	\begin{itemize}
		\item[{\rm (i)}] $\tilde{\mu}_u$ is a completely random measure with intensity
		\begin{equation*}
			\nu^{(u)}(dv,dx) = \mathrm{e}^{-uv} \rho(dv) \alpha(dx) ;
		\end{equation*}
		\item[{\rm (ii)}] for $i=1, \dots, K_n$, the jumps $J_i^{(u)}$ located at the distinct values $X_i^*$ admit density
		\begin{equation*}
			f_i(s) \propto s^{n_i} \mathrm{e}^{-us} \rho(ds) ;
		\end{equation*}
		\item[{\rm (iii)}] for $i=1, \dots, K_n$, $\tilde{\mu}_u$ and $J_i^{(u)}$  are mutually independent.
	\end{itemize}
\end{theorem}

Based on Theorem~\ref{TheoremNRMIposterior}, the posterior characterization of $\tilde P$ is easily deduced and formalized in the next result.

\begin{theorem}[\citealp{MR2508332}, Theorem 2]\label{thm:posterior_p}
	Conditionally on $X^{(n)}$ and $U_n$, $\tilde{P}$ is still an NRMI and it equals in distribution the random probability measure
	\begin{equation*}
		\frac{\tilde{\mu}_{_{U_n}} + \sum_{i=1}^{K_n}J_i^{(U_n)} \delta_{X_i^*}}{\tilde{\mu}_{_{U_n}}(\mathbb{X}) + \sum_{i=1}^{K_n}J_i^{(U_n)}} .
	\end{equation*} 
\end{theorem}

\begin{example}[NGG process] An effective description of the use of Theorems~\ref{TheoremNRMIposterior} and \ref{thm:posterior_p} can be given by taking $\tilde\mu$ to be a normalized generalized gamma process. Consider a generalized gamma CRM as defined in Section 3.2 with $\tau=1$. The normalized generalized gamma (NGG) process is obtained by normalization and its posterior distribution is as follows.

\begin{proposition}[\citealt{MR2730661}]\label{TheoremNGGPosterior}
	Let $U_n$ be a random variable whose density, conditional on $X^{(n)}$, is 
	\begin{equation*}
		f(u \mid X^{(n)}) \propto \frac{u^{n-1} \mathrm{e}^{-\alpha(\mathbb{X}) (1+u)^{\sigma}}}{(u+1)^{n-k\sigma}} .
	\end{equation*}
	The posterior distribution of $\tilde{\mu}$, given $X^{(n)}$, is a mixture with respect to $U_n$ as in \eqref{eq:posterior_general}, where for any $u>0$
	\begin{enumerate}
		\item[{\rm (i)}] $\tilde{\mu}_u$ is a generalized gamma process with intensity
		\begin{equation*}
			\nu^{(u)}(dx,dv) = \frac{\sigma}{\Gamma(1-\sigma)} v^{-1-\sigma} \mathrm{e}^{-(u+1)v} dv \alpha(dx) ;
		\end{equation*}
		\item[{\rm (ii)}] for $i=1, \dots, K_n$, the jumps $J_i^{(u)}$ located at the distinct values $X_i^*$ admit density
		\begin{equation*}
			J_i^{(u)} \sim \textnormal{Ga}(n_i-\sigma, u+1) ;
		\end{equation*}
		\item[{\rm (iii)}] for $i=1, \dots, K_n$, $\tilde{\mu}_u$ and $J_i^{(u)}$ are mutually independent.
	\end{enumerate}
\end{proposition}
The N-IG process, defined in Section~\ref{Section2} via its finite-dimensional distributions, is recovered as a special case for $\sigma=1/2$. Hence, Proposition \ref{TheoremNGGPosterior} immediately yields its posterior characterization.

\end{example}

\section{Construction via predictive distributions}
\label{Section4}

The third construction of the Ferguson--Dirichlet process, which draws inspiration from \cite{MR350949}, relies on the system of predictive distributions. This perspective is closely connected to recent advances on predictive characterizations of nonparametric priors \citep{MR1481784,de2013gibbs}.

\subsection{Predictive characterization}

In the Bayesian nonparametric framework, estimation and prediction are naturally linked: under quadratic loss, the Bayes estimator of $\tilde P$ is its posterior mean, which coincides with the predictive distribution for the next observation. Indeed, by iterated expectation,
\begin{equation} \label{eq:predictive}
	\begin{split} 
		\mathbb{P}(X_{n+1} \in A \mid X^{(n)}) &= \int_{\mathscr{P}} P(A) Q(dP \mid X^{(n)}) \\ 
		& = \mathbb{E}_{Q}(\tilde{P}(A) \mid X^{(n)}) , 
	\end{split} 
\end{equation}
where $Q(\cdot \mid X^{(n)})$ denotes the posterior law of the random probability measure $\tilde P$.

\cite{MR350949} derived the Bayes estimator of the distribution function $\tilde F$ for the Ferguson--Dirichlet process, which takes on a particularly simple and intuitive form. By Proposition \ref{prop:DP_mom}, $\mathbb{E} (\tilde F(t))=P_0((-\infty,t])=F_0(t)$, commonly referred to as the prior guess for $\tilde F$. Let the empirical distribution function of the sample $X^{(n)}$ be
\begin{equation*}
	F_n(t \mid X^{(n)}) = \frac{1}{n} \sum_{i=1}^{n}\delta_{X_i}((-\infty, t]).
\end{equation*}
Then the Bayes estimator of $\tilde F$ is a convex combination of the prior guess and the empirical distribution
\begin{equation}\label{EquationDirichletProcessPredictive1}
	\hat{F}_n	(t \mid X^{(n)}) = \frac{\alpha(\mathbb{R})}{\alpha(\mathbb{R}) + n} \, F_0(t) + \frac{n} {\alpha(\mathbb{R}) + n} \, F_n(t \mid X^{(n)}).
\end{equation}
This follows immediately from the conjugacy property of the Ferguson--Dirichlet process (Theorem  \ref{TheoremDirichletProcessPosterior}). Moreover, by \eqref{eq:predictive}, the predictive distribution satisfies $\mathbb{P}(X_{n+1} \in (-\infty, t] \mid X^{(n)})= \hat{F}_n	(t \mid X^{(n)})$.

Remarkably, predictive distributions that are convex combinations of the prior guess and the empirical measure characterize the Ferguson--Dirichlet process.

\begin{theorem}[\citealt{regazzini1978,MR1130354}]
	Let $(X_n)_{n \geq 1}$ be an exchangeable sequence and assume that, for any $n \geq 1$, the predictive distribution is a convex combination of a probability measure $P_0$ and the empirical measure, that is, for any $A$ in $\mathscr{X}$ and $n\ge 1$
	\begin{equation*}
		\mathbb{P}(X_{n+1} \in A \mid X^{(n)}) = p_n P_0(A) + (1 - p_n)\frac{1}{n}\sum_{i=1}^{n}  \delta_{X_i}(A).
	\end{equation*}
	Then the de~Finetti measure associated with $(X_n)_{n \geq 1}$ is the law of a Ferguson--Dirichlet process with parameter-measure $\alpha$ such that $P_0=\alpha/\alpha(\mathbb{X})$.
\end{theorem}
Note that no assumption is made on the weights $(p_n)_{n\ge1}$, and it is the characterization that implies
\begin{equation*}
	p_n=\alpha(\mathbb{R})/ (\alpha(\mathbb{R}) + n).
\end{equation*}
Thus, within the exchangeable framework, using a convex combination of a prior guess and the empirical measure as the prediction rule is equivalent to adopting the law of a Ferguson--Dirichlet process as the prior distribution.

A definition in terms of predictive distributions is also given in \cite{blackwellmacqueen1973} via a P\'olya urn scheme. They introduce a \emph{P\'olya sequence} $(X_n)_{n\ge 1}$ with parameter $\alpha$, which is generated according to the following scheme
\begin{align*}
	\mathbb{P}(X_1\in A)
	&=\frac{\alpha(A)}{\alpha(\mathbb{X})},\\
	\intertext{and, for any $n\ge 1$,}
	\mathbb{P}(X_{n+1}\in A\,|\,X^{(n)})
	&=\frac{\alpha(A)+\sum_{i=1}^n \delta_{X_i}(A)}{\alpha(\mathbb{X})+n}.
\end{align*}
It is shown that the predictive distribution converges weakly, almost surely, to a limiting random measure $\tilde P^*$, which is a Ferguson--Dirichlet process with parameter measure $\alpha$. Moreover, conditional on $\tilde P^*$, the variables of the P\'olya sequence $(X_n)_{n\ge 1}$ are independent with common distribution $\tilde P^*$. This result was strengthened in two directions by \cite{MR1481784}: (i) the scheme was extended to the broad class of \emph{species sampling models}, which includes the Ferguson--Dirichlet process as a special case; (ii) in this setting the predictive distribution of the generalized P\'olya sequence converges almost surely, in total variation, to the associated species sampling random probability measure.

If $P_0$ is diffuse, the weight $\alpha(\mathbb{R})/ (\alpha(\mathbb{R})+n)$ coincides with $\mathbb{P}(X_{n+1}=\text{``new''} \mid X^{(n)})$, i.e., the probability that the $(n{+}1)$--th observation is a value not previously seen in the sample $X^{(n)}$. This probability depends only on the parameter $\alpha(\mathbb{R})$ and the sample size $n$, and on no other sample information. This implies that $K_n/\log n \rightarrow \alpha(\mathbb{R})$ for $n \rightarrow \infty$. See \cite{MR350950} for further details.

\subsection{Richer predictive structures and Gibbs-type priors}

As in the predictive construction of the Ferguson--Dirichlet process, and more generally for discrete random probability measures with diffuse base measure $P_0$, a key role is played by
\begin{equation}\label{EqProbNewValue}
	\mathbb{P}(X_{n+1}=\text{``new''} \mid X^{(n)}).
\end{equation}
In view of this, \cite{de2013gibbs} characterized predictive distributions of exchangeable sequences according to how the probability of generating a new value depends on the sample information:
\begin{itemize}
	\item[(P1)] The probability \eqref{EqProbNewValue} depends on $n$ but not on the number $K_n$ of distinct values nor on their frequencies $(n_1,\ldots,n_{K_n})$, i.e.
	\begin{equation*}
		\mathbb{P}(X_{n+1}=\text{``new''}\mid X^{(n)}) \;=\; f\bigl(n, \text{model parameters}\bigr).
	\end{equation*}
	\item[(P2)] The probability \eqref{EqProbNewValue} depends on $n$ and $K_n$, but not on $(n_1,\ldots,n_{K_n})$, i.e.
	\begin{equation*}
		\mathbb{P}(X_{n+1}=\text{``new''}\mid X^{(n)}) \;=\; f\bigl(n,\ K_n,\ \text{model parameters}\bigr).
	\end{equation*}
	\item[(P3)] The probability \eqref{EqProbNewValue} depends on all the information in $X^{(n)}$, i.e.
	\begin{equation*}
		\mathbb{P}(X_{n+1}=\text{``new''}\mid X^{(n)}) \;=\; f\bigl(n,\ K_n,\ n_1,\ldots,n_{K_n},\ \text{model parameters}\bigr).
	\end{equation*}
\end{itemize}
Moreover, \cite{de2013gibbs} characterize the members of (P1) and (P2) within species sampling models. (P1) holds if and only if the prior is Ferguson--Dirichlet. By contrast, (P2) holds if and only if the prior is of Gibbs-type, a broad class generalizing the Ferguson--Dirichlet process \citep{MR2160320,lmp2007bmka} and discussed below. For the most general structure (P3), only a few instances are known, such as a specific generalized Ferguson--Dirichlet process \citep{MR1983542,lmp05} and the hierarchical random probability measures considered in \cite{camerlenghi2018}. Unsurprisingly, such generality entails limited analytical tractability. Identifying a class of nonparametric priors in (P3) that satisfies Ferguson's two \emph{desiderata} remains an open problem.

For the Ferguson--Dirichlet process, \eqref{EqProbNewValue} depends solely on the sample size, which entails restrictive features such as the logarithmic growth of $K_n$. This can be limiting in several application domains, including linguistics and topic modeling \citep{teh2006,teh2010}, networks \citep{caron2017,caron2020}, population genetics \citep{feng2010,ruggiero2017}, and species sampling problems \citep{lmp2007bmka,flp2009}. This motivates Gibbs-type priors, corresponding to class (P2) above, which represent an important and, as argued in \cite{de2013gibbs}, perhaps the most natural generalization of the Ferguson--Dirichlet process. Like the Ferguson--Dirichlet process, Gibbs-type priors can also be characterized by their predictive structure, a result first stated in \cite{lmp2007bmka} and derived from \cite{MR2160320}.

\begin{proposition}\label{thm:gibbs_characterize}
	Let $\tilde P=\sum \tilde p_i \delta_{V_i}$ be a discrete random probability measure with probability weights $(\tilde p_i)_{i \geq 1}$ independent of the atoms $(V_i)_{i \geq 1}$, which are i.i.d. from a diffuse $P_0$. The law of $\tilde P$ is a Gibbs-type prior with parameters $\sigma \in (-\infty, 1)$ and $P_0$ if and only if it induces predictive distributions of the form 
	\begin{equation*}
		\mathbb{P}(X_{n+1} \in A \mid X^{(n)})
			= \frac{V_{n+1, K_n+1}}{V_{n, K_n}} P_0(A) + \left( 1 - \frac{V_{n+1, K_n+1}}{V_{n, K_n}} \right) \frac{\sum_{i=1}^{K_n}(n_i - \sigma)\delta_{X_i^*}(A)}{n - \sigma K_n} ,
	\end{equation*}
	for any $n\ge 1$ and $A$ in $\mathscr{X}$, where the weights $\{ V_{n,k}:\:  n\geq 1, \, 1\leq k\leq n \}$ satisfy the following recursion
	\begin{equation*}
		V_{n,k} = (n - k\sigma)V_{n+1, k} + V_{n+1, k+1} .
	\end{equation*}
\end{proposition}
Therefore, a Gibbs-type prior is characterized by three components: the prior guess $P_0=\mathbb{E}(\tilde P)$, the parameter $\sigma<1$ and the weights $V_{n,k}$. Moreover, one has
\begin{equation*}
	\mathbb{P}(X_{n+1}=\text{``new''} \mid X^{(n)}) = \frac{V_{n+1, K_n+1}}{V_{n, K_n}},
\end{equation*}
which depends on both the sample size $n$ and the number of distinct values in the sample $K_n$, as required by (P2). This predictive structure is clearly more flexible than that of the Ferguson--Dirichlet process; for example, when $\sigma\in(0,1)$ it yields $K_n$ growing at rate $n^{\sigma}$. See \cite{pitman2006,de2013gibbs} for details.

As for Ferguson's first desideratum on large support, Gibbs-type priors satisfy it. In particular, as shown in \cite{dlp2013}, the support of a Gibbs-type prior is the set of all probability measures whose support is contained in that of $P_0$, when either $\sigma \ge 0$ or $\sigma < 0$ and a sufficient condition is met. Regarding analytical tractability, this desideratum is only partially met for general Gibbs-type priors, because solving the recursion in Proposition~\ref{thm:gibbs_characterize} to obtain fully explicit prediction rules is challenging. Nevertheless, several important special cases enjoy full analytical tractability.

The most popular example is the Pitman-Yor process \citep{MR1434129}, characterized by $P_0$ and parameters $(\theta,\sigma)$ with the admissible ranges: if $\sigma \in [0,1)$ then $\theta > -\sigma$, while if $\sigma < 0$ one has $\theta = r |\sigma|$ for some positive integer $r$. The corresponding weights that satisfy the recursion in Proposition~\ref{thm:gibbs_characterize} are
\begin{equation*}
	V_{n,k} = \frac{\prod_{i=1}^{k-1}(\theta + i\sigma)}{(\theta+1)_{n-1}},
\end{equation*}
and they yield the simple and intuitive predictive rule
\begin{equation*}
	\mathbb{P}(X_{n+1} \in A \mid X^{(n)}) = \frac{\theta+\sigma K_n}{\theta+n} P_0(A) + \frac{n-\sigma K_n}{\theta+n} \ \sum_{i=1}^{K_n}\frac{n_i - \sigma}{n - \sigma K_n} \ \delta_{X_i^*}(A)
\end{equation*}
for any $n \ge 1$ and $A \in \mathscr{X}$. Clearly, when $\sigma = 0$ these expressions reduce to the predictive distributions of the Ferguson--Dirichlet process. The directing de~Finetti measure is then a Pitman-Yor process prior with parameters $(\sigma,\theta)$. Although the Pitman-Yor process is not an NRMI, it admits both \emph{a priori} and \emph{a posteriori} representations in terms of transformed CRMs \citep{MR1434129,pitman2006,MR2730661}. Consequently, approximate simulation of its realizations is feasible by sampling its largest atom weights \citep{MR4216518,MR4764755}.

Another notable example of a Gibbs-type prior is the NGG process, introduced in \cite{MR2370077} and already discussed as a special case of an NRMI in Section~\ref{sec:NRMI}. In fact, the only nonparametric prior being simultaneously an NRMI and of Gibbs--type is the NGG process \citep{lpw2008}. The NGG corresponds to a Gibbs--type prior with weights
\begin{equation*}
	V_{n,k} = \frac{\mathrm{e}^{\beta}\sigma^{k-1}}{\Gamma(n)} \sum_{i=0}^{n-1}\binom{n-1}{i}(-1)^i\beta^{i/\sigma}\Gamma(k-i/\sigma;\beta)
\end{equation*}
with $\beta>0$, $\sigma\in(0, 1)$ and $\Gamma(x,a)$ denoting the incomplete gamma function, which also satisfy the recursion in Proposition~\ref{thm:gibbs_characterize}. In particular, if $\sigma=1/2$, we
recover the normalized inverse Gaussian process presented in Section~\ref{sec:NIG}. It is worth noting that the NGG weights were not obtained by solving the recursion in Proposition~\ref{thm:gibbs_characterize}. Instead, they were derived by starting from an exchangeable sequence directed by an NGG prior, then computing the induced exchangeable random partition and the resulting predictive distribution. This indirect route appears more promising than attempting to solve the recursion directly to obtain further explicit Gibbs--type special cases. Identifying additional instances remains an interesting open problem.

\begin{rmk}
	The three constructions of the Ferguson--Dirichlet process reviewed in Sections~\ref{Section2}, \ref{Section3}, and~\ref{Section4} are equivalent but lead to rather different generalization paths. Each construction highlights a different aspect of the process and comes with its own advantages and limitations. The finite-dimensional construction in Section~\ref{Section2} is closest to Ferguson's original definition and makes the role of additivity transparent. Its main limitation is that additivity is a very restrictive property: besides the Dirichlet case, only a few explicit additive distributions on the simplex are known, as discussed in Section~\ref{Section2}. The normalization construction in Section~\ref{Section3}, instead, is much more flexible. By replacing the gamma process with a general CRM, one obtains the broad class of NRMIs, for which posterior characterizations and simulation methods can be studied through the L\'evy intensity and the Laplace functional. Finally, the predictive construction in Section~\ref{Section4} is especially useful for understanding the learning mechanism induced by the prior, since it directly displays the reinforcement structure and the probability of observing new values. This perspective naturally leads to Gibbs--type priors, although directly identifying explicit prediction rules beyond the Ferguson--Dirichlet process is technically challenging because of the coherence requirements imposed by exchangeability; see, for instance, \cite{fortini}. Thus, while the three constructions coincide for the Ferguson--Dirichlet process, the generalizations they suggest differ substantially in scope, tractability, and statistical interpretation.
\end{rmk}

\section{Dependent Ferguson--Dirichlet processes}
\label{Section5}

From \cite{MR350949} through developments up to the late 1990s, Bayesian nonparametric applications were largely confined to settings in which the observations were assumed exchangeable. However, exchangeability is not appropriate when data exhibit heterogeneity, for example when covariates affect a response. This occurs in topic modeling, meta-analysis, two-sample problems, nonparametric regression, and inference with time-dependent data or change-points. The limitations of exchangeability were already noted in \cite{de1938condition} and are effectively described in the following passage:
\begin{quote}
	\textsl{But the case of exchangeability can only be considered as a limiting case: the case in which this ``analogy'' is, in a certain sense, absolute for all events under consideration. To get from the case of exchangeability to other cases which are more general but still tractable, we must take up the case where we still encounter ``analogies'' among the events under consideration, but without attaining the limiting case of exchangeability.}
\end{quote}
In applications with covariate-dependent data, a form of probabilistic symmetry more general than exchangeability is required. Let $z\in\mathcal{Z}$ be a covariate, and introduce a collection of (possibly dependent) random probability measures $\bm{P}_{\mathcal{Z}}=\{\tilde P_z:\: z\in\mathcal{Z}\}$ assuming 
\begin{equation}\label{eq:partial_exchange}
	X_{z_1,i_1},\ldots,X_{z_k,i_k}|\bm{P}_{\mathcal{Z}}
	\sim\tilde P_{z_1}\,\times\,\cdots\,\times\,\tilde P_{z_k}
\end{equation}
for any $k\ge 1$, $z_1,\ldots,z_k\in\mathcal{Z}$ and $i_1,\ldots,i_k\ge 1$. A simple, yet important, case is when $\mathcal{Z}=\{1,\ldots,d\}$ since this can be seen as a situation where data are recorded under different, though related, experimental conditions. Hence, the above dependence structure is equivalent to assuming
\begin{itemize}
	\item[(C1)] Homogeneity within each experimental condition;
	\item[(C2)] Heterogeneity across different experimental conditions.
\end{itemize}
When $\mathcal{Z}$ is finite, the conditional independence in \eqref{eq:partial_exchange} is equivalent to \textit{partial exchangeability} \citep{de1938condition}. We give the definition for $d=2$, the extension to general $d$ being immediate.
\begin{definition}[\citealp{de1938condition}]
	The array $(\textbf{X}_1, \textbf{X}_2) = (X_{1,i} ,X_{2,j})_{ i,j\geq 1 }$ is \emph{partially exchangeable} if
	\begin{equation*}
		(X_{1,1}, \ldots, X_{1,n_1}, X_{2,1}, \ldots, X_{2,n_2}) \eqd (X_{1,\pi(1)}, \ldots, X_{1,\pi(n_1)}, X_{2,\phi(1)}, \ldots, X_{2,\phi(n_2)}) ,
	\end{equation*}
	for any $n_1, n_2 \geq 1$ and any permutations $\pi$ and $\phi$ of $(1, \ldots, n_1)$ and $(1, \ldots, n_2)$, respectively.
\end{definition}
Hence, the distribution of a partially exchangeable array is invariant with respect to permutations within each condition, but not necessarily to permutations of elements associated with different conditions. This corresponds to the mathematical formalization of (C1) and (C2) above.  As in the exchangeable case, a representation theorem holds.
\begin{theorem}[\citealp{de1938condition}]
	The array $(\textbf{X}_1, \textbf{X}_2)$ is partially exchangeable if and only if there exists some probability measure $Q$ on $\mathcal{P}^2$ such that
	\begin{equation}\label{eq:represent_part_exchange}
		\mathbb{P}\Big(X_1^{(n_1)}\in \mathsf{X}_{i=1}^{n_1}A_i, \:
			X_2^{(n_2)}\in \mathsf{X}_{j=1}^{n_2} B_j \Big) = \int_{\mathcal{P}^2} \prod_{i=1}^{n_1} P_1(A_i) \prod_{j=1}^{n_2} P_2(B_j) Q(d P_1, d P_2) 
	\end{equation}
	for any positive integers $n_1$ and $n_2$, and sets $A_1,\ldots,A_{n_1}$ and $B_{1},\ldots,B_{n_2}$ in $\mathscr{X}$, where
	$X_i^{(n_i)}=(X_{i,1},\ldots,X_{i,n_i})$ for $i=1,2$.
\end{theorem}
Equivalently, the mixture representation \eqref{eq:represent_part_exchange} can be expressed in terms of a generative model of the form \eqref{eq:partial_exchange} with $\mathcal{Z}=\{1,2\}$ and
\begin{equation*}
		(\tilde{P}_1, \tilde{P}_2)	\sim  Q ,
\end{equation*}
where $(\tilde{P}_1, \tilde{P}_2)$ is a vector of dependent random probability measures and $Q$ is the de~Finetti measure acting as a prior. One should still aim for a $Q$ that achieves Ferguson's two desiderata, namely large support and analytical tractability. However, in this more general framework the dependence between $\tilde P_1$ and $\tilde P_2$ induced by $Q$  is itself central for inference, which motivates an additional
\begin{quote}
	\textit{Desideratum.} Specify a $Q$ ensuring a wide range of dependence structures between $\tilde P_1$ and $\tilde P_2$.
\end{quote}
Pursuing this objective calls for principled measures of dependence between $\tilde P_1$ and $\tilde P_2$, so as to understand the dependence induced by $Q$. This has motivated a recent research stream aiming at introducing measures of dependence relying on optimal transport theory \citep{marta2021,marta_jasa} and reproducing kernel Hilbert spaces \citep{marta_arxiv}; see also \cite{catalanolavenant2025} for a review.

An early extension of the Ferguson--Dirichlet process to this setting is the mixture of products of Ferguson--Dirichlet processes in \cite{cifreg78}, which is tailored to the case of a finite $\mathcal{Z}$ space. The most influential approach, applicable to general covariate spaces $\mathcal{Z}$, relies on dependent stick-breaking constructions for the Ferguson--Dirichlet process. This rich line of research originates in \cite{maceachern1999dependent} and in the 2000 manuscript now published as \cite{10.1214/26-BA1610}, where the dependent Ferguson--Dirichlet process was introduced. See \cite{mueller2022} and \cite{wade_inacio} for recent stimulating reviews. In line with the CRM construction of the Ferguson--Dirichlet process presented in Section \ref{Section3}, we close by pointing to random-measure based approaches for defining dependent Ferguson--Dirichlet processes and more general classes of dependent processes, organized by their probabilistic structure:
\begin{itemize}
	\item[(i)] \emph{Additive structures.} This approach was originally introduced in \citep{MR2088779} for the Ferguson--Dirichlet case. The underlying theory has been investigated in \cite{MR3217444, MR3131980} in the more general framework of NRMIs.
	\item[(ii)] \emph{Hierarchical structures.} This approach was pioneered by \cite{MR2279480}, who introduced a hugely successful version of the dependent Ferguson--Dirichlet process, known as the \emph{hierarchical Dirichlet process}. The distribution theory for this class of priors has been developed, and extended beyond the Ferguson--Dirichlet case, in \cite{MR3909927,MR4255111,MR4827151,10.1214/26-BA1609}.
	\item[(iii)] \emph{Nested structures.} The first proposal of a dependent Ferguson--Dirichlet process with a nested structure can be found in \cite{MR2528831} and it is known as the \emph{nested Dirichlet process}. Variants and extensions of this model, aimed at overcoming some of its shortcomings, can be found in \cite{MR4044854}, \cite{MR4558734} and \cite{denti23}.
\end{itemize} 

To sum up, most constructions of dependent nonparametric priors build either on transformations of random measures, as in \cite{MR350949}, or on stick-breaking constructions. Stick-breaking is especially convenient computationally, while the random-measure viewpoint, originating in Ferguson's work, often facilitates \emph{a priori} and \emph{a posteriori} distributional analysis. We are not aware of dependent versions of the Ferguson--Dirichlet process (or possible generalizations) defined via finite-dimensional distributions as in Section~\ref{Section2} or via prediction rules as in Section~\ref{Section4}.

\section{Further perspectives and concluding remarks}\label{SectionRemarks}

Over the past five decades, theoretical and methodological advances in Bayesian Nonparametrics have been profoundly shaped by Ferguson's 1973 breakthrough paper introducing the Ferguson--Dirichlet process. Ferguson's formulation provided a rigorous probabilistic foundation for the use of random measures in modern Bayesian nonparametric inference. The many generalizations it inspired are not only alternatives to the Ferguson--Dirichlet process, but also help identify which features of Ferguson's construction are exceptional and which ones can be preserved under greater modeling flexibility. These developments have been instrumental in making Bayesian nonparametric modeling popular across Statistics and Machine Learning, with the Ferguson--Dirichlet process serving as the key reference model. Comprehensive accounts of developments following \cite{MR350949} can be found in several excellent and inspiring monographs, including \cite{ghosh_rama},  \cite{HHMW2010}, \cite{mueller2015} and \cite{ghosal2017}. The motivation for considering these alternatives is not merely to enlarge the catalogue of Bayesian nonparametric priors. In species sampling and clustering, alternatives to the Ferguson--Dirichlet process allow richer reinforcement mechanisms and different growth rates for the number of distinct values. In mixture models, the prior on the mixing distribution affects the induced partition of the latent variables and hence the number and relative sizes of the mixture components. These examples illustrate the basic trade-off: the Ferguson--Dirichlet process remains exceptional in its simplicity and tractability, while its generalizations are useful when additional predictive or structural flexibility is needed.

In the previous sections, we reviewed the constructions in \cite{MR350949}, which have inspired numerous developments beyond the original Ferguson--Dirichlet framework. Other equivalent formulations have also appeared in the literature.

The most prominent is the stick-breaking construction of \cite{sethuraman94}, which represents a Ferguson--Dirichlet process as
\begin{equation*}
	\tilde P=\sum_{i\ge1}\tilde p_i\delta_{V_i},
	\qquad
	\tilde p_i=W_i\prod_{j<i}(1-W_j),
\end{equation*}
where $W_i\stackrel{\mathrm{ind}}{\sim}\mathrm{Beta}(1,\alpha(\mathbb X))$ and $V_i\stackrel{\mathrm{iid}}{\sim}P_0=\alpha/\alpha(\mathbb X)$, independently of the weights. This representation made the process especially accessible by expressing it through independent beta breaks and i.i.d. atoms, thereby enabling direct simulation, finite truncations, and dependent stick-breaking extensions. Stick-breaking representations have also been developed for other priors discussed in the previous sections, including the Pitman--Yor process \citep{PPY92}, the normalized inverse Gaussian process \citep{FLP12}, and homogeneous NRMIs \citep{Favaro2016}. More broadly, stick-breaking ideas have generated a rich literature on alternative and dependent constructions; see, e.g. \cite{MR1481784,ishwaran2001,dunson2008,gil2023,li2025,10.1214/25-AOS2607}.

Another notable formulation, due to \cite{MR438568}, represents the Ferguson--Dirichlet process as a neutral-to-the-right (NTR) prior \citep{doksum1974}.  Using the posterior characterization of NTR priors \citep{ferguson1979}, one can show that the Ferguson--Dirichlet process is not conjugate in the presence of censored observations. This lack of conjugacy was also established, outside the NTR framework, by \cite{susarla76}. Moreover, since NTR processes can be recast as priors for cumulative hazards \citep{ramam2003}, the Ferguson--Dirichlet process is also recovered as a special case of the prior induced by a beta process on the cumulative hazard \citep{hjort1990}. 

In concluding this exposition of the fundamental influence of the Ferguson--Dirichlet process on the evolution of Bayesian Nonparametrics, it is fitting to recall a passage from \cite{MR438568} that, in retrospect, remains strikingly relevant today, even though the MCMC revolution was yet to come.
\begin{quote}
	\textsl{One of the drawbacks of decision theory in general and of the Bayesian approach to it in particular is the difficulty of putting the cost of the computation into the model. This drawback is particularly severe in Bayesian nonparametric problems. There are no doubt examples in which ``quick and easy'' rules are preferable to ``optimal'' rules for a Bayesian simply because it costs less to perform the computations. On the other hand, Bayes rules are certainly desirable since generally they are admissible and have nice large sample properties. Therefore, it behooves the statisticians to suggest large classes of easily computable Bayes rules in the hope that users may find some rules to their liking.}
\end{quote}
The enduring legacy of \cite{MR350949} is not only the Ferguson--Dirichlet process but a rigorous, principled foundation for Bayesian nonparametric modeling, with clear desiderata, constructions that meet them, and an explicit posterior-predictive calculus. Since the advent of MCMC, this framework has been complemented by a sustained search for easily computable Bayes rules, with analytical tractability becoming essential to efficient algorithms. This combination has made Bayesian nonparametric models viable in complex applied settings. The path ahead is to keep this balance, pursuing richer predictive schemes, broader classes with large support, and scalable algorithms, while recognizing the process that still anchors both theory and practice: the Ferguson--Dirichlet process.

\begin{appendix}
\section{Outline of the proof of Proposition \ref{TheoremAprioriMoments}}

A key step for the derivation of  $\mathbb{E}[\tilde P(A)]$ consists in rewriting the normalization as
\begin{equation*}
	\tilde{P}(A) = \frac{\tilde{\mu}(A)}{\tilde{\mu}(\mathbb{X})} = \int_{0}^{\infty} \tilde{\mu}(A) \mathrm{e}^{-u \tilde{\mu}(\mathbb{X})} du .
\end{equation*}
Then we have
\begin{equation*}
	\begin{split}
		\mathbb{E}(\tilde{P}(A)) &= \mathbb{E} \left( \int_{0}^{\infty} \tilde{\mu}(A) \mathrm{e}^{-u \tilde{\mu}(\mathbb{X})} du \right) \\
		&= \mathbb{E} \left( \int_{0}^{\infty} \tilde{\mu}(A) \mathrm{e}^{-u \tilde{\mu}(A)} \mathrm{e}^{-u \tilde{\mu}(A^C)} du \right) \\
		&= \int_{0}^{\infty} \mathbb{E} \left( \tilde{\mu}(A) \mathrm{e}^{-u \tilde{\mu}(A)} \right) \mathbb{E}\left( \mathrm{e}^{-u \tilde{\mu}(A^C)} \right) du ,
	\end{split}
\end{equation*}
where the last equality follows from Fubini's theorem and the independence between $\tilde{\mu}(A)$ and $\tilde{\mu}(A^C)$. Also, using the Laplace transform of $\tilde{\mu}$, we have
\begin{equation*}
	\mathbb{E}\left( \mathrm{e}^{-u \tilde{\mu}(A^C)} \right) = \exp \left( -\int_{0}^{\infty}(1 - \mathrm{e}^{-uv})\rho(dv)\alpha(A^C) \right) .
\end{equation*}
Denote by 
\begin{equation}
	\psi(u) := \int_{0}^{\infty}(1 - \mathrm{e}^{-uv})\rho(dv) .
\end{equation}
It follows that
\begingroup
\allowdisplaybreaks
\begin{align*}
		\mathbb{E}(\tilde{P}(A))
		&= \int_{0}^{\infty} \mathbb{E} \left( -\frac{d}{du} \mathrm{e}^{-u \tilde{\mu}(A)} \right) \mathrm{e}^{-\psi(u)\alpha(A^C)} du \\
		&= \int_{0}^{\infty} \left( -\frac{d}{du} \mathbb{E} (\mathrm{e}^{-u \tilde{\mu}(A)}) \right) \mathrm{e}^{-\psi(u)\alpha(A^C)} du \\
		&= \int_{0}^{\infty} \left( -\frac{d}{du} \mathrm{e}^{-\psi(u)\alpha(A)} \right) \mathrm{e}^{-\psi(u)\alpha(A^C)} du \\
		&= \alpha(A) \int_{0}^{\infty} \left( \frac{d}{du} \psi(u) \right) \mathrm{e}^{-\psi(u)\alpha(A)} \mathrm{e}^{-\psi(u)\alpha(A^C)} du \\
		&= \alpha(A) \int_{0}^{\infty} \left( \frac{d}{du} \psi(u) \right) \mathrm{e}^{-\psi(u)\alpha(\mathbb{X})} du \\
		&= \frac{\alpha(A)}{\alpha(\mathbb{X})} \int_{0}^{\infty}  -\frac{d}{du} \mathrm{e}^{-\psi(u)\alpha(\mathbb{X})} du \\
		&= \frac{\alpha(A)}{\alpha(\mathbb{X})} (1 - \mathrm{e}^{-\rho(\mathbb{R}^+)\alpha(\mathbb{X})}) \\
		&= \frac{\alpha(A)}{\alpha(\mathbb{X})},
\end{align*}
\endgroup
where the last equality holds by condition (H2) in Section~\ref{sec:NRMI}, which ensures that $\tilde\mu(\mathbb{X})>0$ almost surely. See \cite{MR2255112} for a detailed justification of the displayed equalities. The variance and covariance can be derived in a similar way. We refer the reader to Appendix A of \cite{MR2255112}.

\end{appendix}

\section*{Acknowledgments}
This paper grew out of the talk \emph{``Ferguson's 1973 Breakthrough: A Timeless Inspiration''} given by I.~Pr\"unster at the event \emph{``50 Years of Tom Ferguson's 1973 Dirichlet Process Paper \& 25th Anniversary Celebration''} held at the Department of Statistics, UCLA (November 2023). It was a great honor to have Prof.~Thomas S. Ferguson present at the event. We also wish to thank the participants for their valuable comments and express our warm gratitude to Michele Guindani and Mark Handcock for their outstanding efforts in organizing the celebration. We are also grateful to Eugenio Regazzini for drawing our attention to \cite{definetti1935}. Part of this work was carried out while J. Zhang was affiliated with Bocconi University. All authors contributed equally to this work.

\section*{Funding}
The authors were supported by the European Union-Next Generation EU PRIN-PNRR (Project P2022H5WZ9). 

\bibliographystyle{abbrvnat}
\bibliography{bibliography.bib}

\end{document}